\documentclass[11pt,reqno]{amsart}
\usepackage[T1]{fontenc}
\usepackage{amsmath,amssymb,amsthm,mathtools,esint}
\usepackage[a4paper,margin=28mm]{geometry}
\usepackage{needspace}
\usepackage[nopatch=footnote]{microtype}
\usepackage[colorlinks=true,linkcolor=blue,citecolor=blue,urlcolor=blue]{hyperref}
\AtBeginDocument{%
  \setlength{\parskip}{2pt plus .5pt minus .5pt}%
  \setlength{\jot}{4pt}%
}
\newtheorem{theorem}{Theorem}[section]
\newtheorem{proposition}[theorem]{Proposition}
\newtheorem{lemma}[theorem]{Lemma}

\theoremstyle{remark}

\numberwithin{equation}{section}
\newcommand{\norm}[1]{\left\lVert#1\right\rVert}
\newcommand{\avg}{\fint}
\newcommand{\Q}{\mathcal Q}
\newcommand{\Pcal}{\mathcal P}
\newcommand{\R}{\mathbb R}
\newcommand{\D}{\mathcal D}
\DeclareMathOperator{\Dom}{Dom}
\DeclareMathOperator{\Cov}{Cov}

\title[Heat observability from measurable sets]{Observability from Measurable Sets for Heat Equations with Bounded Space-Time Potentials}

\author[M. Wang]{Ming Wang}
\address[M. Wang]{School of Mathematics and Statistics, HNP-LAMA,
Central South University, Changsha, Hunan 410083, P.R. China}
\email{m.wang@csu.edu.cn}
\author[C. Zhang]{Can Zhang}
\address[C. Zhang]{School of Mathematics and Statistics,
Wuhan University, Wuhan 430072, P.R. China}
\email{canzhang@whu.edu.cn}
\hypersetup{
  pdfauthor={Ming Wang and Can Zhang},
  pdftitle={Observability for Heat Equations with Bounded Space-Time Potentials on Bounded Domains and the Whole Space}
}
\date{November 12, 2026}

\begin{document}

\begin{abstract}
We establish final-state observability for heat equations with real-valued
bounded space--time potentials on bounded \(C^2\) domains, from arbitrary
measurable space--time sets of positive measure. On \(\R^n\), we prove
observability from time-dependent sets whose spatial slices are uniformly
thick on a set of times of positive measure. A counterexample for the free
heat equation on \(\R\) shows that integral thickness alone is insufficient
for observability.
\end{abstract}
\maketitle

\section{Introduction}\label{sec-intro}

\subsection{Main results}

Let $\Omega\subset\R^n$, $n\ge1$, be a bounded connected domain with $C^2$
boundary, and let $T>0$. We consider
\begin{equation}\label{equ-1031-1}
 \begin{cases}
  \partial_tu-\Delta u+a(x,t)u=0&\text{in }\Omega\times(0,T),\\
  u=0&\text{on }\partial\Omega\times(0,T),\\
  u(\cdot,0)=u_0&\text{in }\Omega,
 \end{cases}
\end{equation}
where $a$ is real valued and
\begin{equation}\label{equ-1031-2}
 a\in L^\infty(\Omega\times(0,T)),\qquad \norm a_\infty\le M.
\end{equation}
For initial data in $L^2(\Omega)$, the solution lies in the energy class
\[
 C([0,T];L^2(\Omega))\cap L^2(0,T;H_0^1(\Omega)).
\]

\begin{theorem}\label{thm-ob-bound}
For every measurable $D\subset\Omega\times(0,T)$ with $|D|\ge d>0$,
every solution of \eqref{equ-1031-1}--\eqref{equ-1031-2} satisfies
\begin{equation}\label{equ-1031-3}
 \norm{u(T)}_{L^2(\Omega)}
       \le C\int_D|u(x,t)|\,\mathrm{d}x\,\mathrm{d}t.
\end{equation}
Here $C=C(\Omega,T,M,d)>0$.
\end{theorem}

Observability from arbitrary positive-measure space--time sets was proved
for the free heat equation in \cite{AEWZ} and for parabolic equations with
time-independent analytic coefficients in \cite{EMZ}. For the free heat
equation, spatial observation sets may also have zero Lebesgue measure
when observations are measured in a supremum norm, under conditions on
Hausdorff dimension \cite{BM,GLMO} or finer log-type Hausdorff content
\cite{HWW}. For autonomous rough
operators, spectral methods give measurable-set observation estimates for
Lipschitz metrics \cite{BM} and real bounded spatial potentials \cite{LBM}.
For bounded space--time potentials, earlier results use open spatial
observations over measurable time sets \cite{PW,PWZ}; refined costs for
Lipschitz potentials and one-dimensional autonomous product observations
are studied in \cite{ZZ}. The case of arbitrary positive-measure space--time
sets was recently established on a bounded interval for real bounded
nonautonomous potentials \cite{XZ}. Theorem~\ref{thm-ob-bound} extends this
result to every spatial dimension, uniformly under a prescribed potential
bound and a positive lower bound for the observation measure.

For the equation in the whole space, a measurable set $\omega\subset\R^n$ is
called $(\gamma,\ell)$-thick if
\begin{equation}\label{equ-1031-4}
 |\omega\cap(x+[0,\ell]^n)|\ge\gamma\ell^n
 \quad\text{for every }x\in\R^n,
 \qquad 0<\gamma\le1,\quad\ell>0.
\end{equation}
The observation set may depend on time, with these parameters uniform
on a common set of observation times.

For a measurable space--time set $D$, we write $D_t=\{x:(x,t)\in D\}$.

\begin{theorem}[Whole-space observability]\label{thm-ob-Rn}
Let $a\in L^\infty(\R^n\times(0,T))$ be real with $\norm a_\infty\le M$,
and let $D\subset\R^n\times(0,T)$ be measurable.
If $D_t$ is $(\gamma,\ell)$-thick for almost every $t$ in a measurable
$\mathcal A\subset(0,T)$ with $|\mathcal A|>0$, then every energy solution of
\begin{equation}\label{equ-1031-5}
 \partial_tu-\Delta u+a(x,t)u=0\quad\text{in }\R^n\times(0,T),\qquad
 u(0)=u_0\in L^2(\R^n),
\end{equation}
satisfies
\begin{equation}\label{equ-1031-6}
 \norm{u(T)}_2\le C\int_{\mathcal A}\norm{u(t)}_{L^2(D_t)}\,\mathrm{d}t,
\end{equation}
where $C=C(n,T,M,\gamma,\ell,|\mathcal A|)$.
\end{theorem}

Thickness characterizes fixed observable sets for the free heat equation
on $\R^n$ \cite{EV,WWZZ}, and measurable thick sets also give observability
for bounded spatial potentials \cite{LBM}. For spatially analytic
space--time potentials, thick-set observability is also known for
fractional parabolic equations \cite{WZ23}. For general bounded space--time
potentials, previous results use spatial regions containing uniformly
distributed balls \cite{DWZ}. Theorem~\ref{thm-ob-Rn} allows arbitrary
measurable thick slices that may depend on time.

Theorem~\ref{thm-ob-bound} and Theorem~\ref{thm-ob-Rn} naturally raise
the question whether integral thickness, which averages spatial density
in time, is necessary and sufficient for space--time observability of the
free heat equation on $\R^n$. The sufficiency question for exact
null-controllability from moving supports was raised in \cite{BEP} and
explicitly recalled as open in \cite[Remark 2.3]{AM}.
It is shown that uniformly thick moving slices give
observability \cite{BEP}, whereas integral thickness characterizes
approximate null controllability with uniform cost \cite{AM}.
Theorem~\ref{thm-integral-counterexample} shows that integral thickness does
not suffice for exact observability, even on $\R$.

\subsection{Main difficulties and proof strategy}

Our interest in measurable observation sets comes from control theory and
from the geometry of observable sets. In optimal control, the bang--bang
property means that an optimal control saturates its admissible bound
almost everywhere. Observability on sets of positive measure is used in
the duality and perturbation arguments proving this property, in
particular for time-optimal heat control \cite{AEWZ,PWZ}. From a geometric
viewpoint, we seek necessary and sufficient conditions under which a set
observes every final state. The thickness characterization on $\R^n$ and
the difference between fixed and moving observations show the importance
of allowing observation sets without interior.

The main difficulty comes from the low regularity and time dependence of
the potential. For equations with analytic coefficients, quantitative analyticity and
propagation of smallness for analytic functions allow one to pass from a
measurable subset to a whole ball \cite{AEWZ,EMZ}. For a merely bounded
space--time potential, this analytic structure is no longer available.

For potentials depending only on space, a different approach is available.
Doubling-index estimates give propagation of smallness for gradients on
sets of sufficiently large Hausdorff dimension \cite{LM}. By applying
these estimates to an elliptic extension, one obtains spectral
inequalities for a fixed spatial operator \cite{BM,LBM}. These inequalities,
together with high-frequency heat decay, give observability \cite{LR}.
For a bounded nonautonomous potential, however, no fixed spatial operator
generates the evolution, and the resulting coupling of modes is not
directly controlled by these spectral estimates.

We outline the proof of Theorem~\ref{thm-ob-bound}.

We use local Hermite estimates \cite{WZ} and averaging over Gaussian
centers to reduce the past frequency on a cube that retains positive
observation density. Gaussian convolution in the center variable then
gives terminal cylinder doubling: for centers $x$ in the cube and
sufficiently small $r$ depending on $\Lambda_{s,t}$,
\begin{equation}\label{equ-1110-1}
 \avg_{\Q_{2r}(x,t)}|u|^2\le C\avg_{\Q_r(x,t)}|u|^2,
\end{equation}
with $C$ independent of $\Lambda_{s,t}$. Caloric polynomial approximation
and the Remez inequality give propagation of smallness on the cube.
Together with the small-ball mass bound, this yields, for every measurable
$E\subset\Omega$ with $|E|\ge m>0$,
\begin{equation}\label{equ-1110-2}
 \norm{u(t)}_2
 \le C_m\exp\!\left(C_m\Lambda_{s,t}^{3/2}\right)
      \norm{u(t)}_{L^2(E)}.
\end{equation}
Applying this estimate to a low-amplitude subset of each slice $D_t$ with
$|D_t|\ge d/(2T)$ gives an $L^1(D_t)$ observation bound. Since $3/2<2$,
the weighted telescoping argument in Proposition~\ref{prop-subquadratic}
then gives
\[
 \norm{u(T)}_2\le C\int_D|u(x,t)|\,\mathrm{d}x\,\mathrm{d}t.
\]

For the proof on $\R^n$, we first select cubes with controlled local energy
ratios which retain a fixed fraction of the terminal energy. By summing the
local estimates and applying the same time argument, we obtain the
whole-space result.


\subsection{Notation}
We write $u(t)=u(\cdot,t)$ and $\norm{\cdot}_2$ for the spatial $L^2$ norm.
 We use Euclidean balls $B_r(x)$, cubes
$Q_r(x)=x+(-r,r)^n$, cylinders
\[
 \Pcal_r(x,t)=B_r(x)\times(t-r^2,t],\qquad
 \Q_r(x,t)=Q_r(x)\times(t-r^2,t].
\]
We also use the averages $\avg_A v=|A|^{-1}\int_Av$ for $0<|A|<\infty$. We set
\[
 G_\sigma(z)=(4\pi\sigma)^{-n/2}\mathrm{e}^{-|z|^2/(4\sigma)},\qquad
 \mathrm{d}\gamma(y)=G_1(y)\,\mathrm{d}y.
\]
We denote the $L^2(\mathrm{d}\gamma)$ inner product and norm by
$\langle\cdot,\cdot\rangle_\gamma$ and $\norm{\cdot}_\gamma$.
All logarithms are natural; $X\lesssim Y$ means $X\le CY$, and
$X\approx Y$ means $C^{-1}Y\le X\le CY$.

For \(0\le s<t\le T\) with \(\norm{u(t)}_2>0\), set
\begin{equation}\label{equ-1031-7}
 \delta=t-s,\qquad
 \Theta_{s,t}=\frac{\mathrm{e}^{2M(t-s)}\norm{u(s)}_2^2}{\norm{u(t)}_2^2},
 \qquad \Lambda_{s,t}=1+\delta^{-1}+\log\Theta_{s,t}.
\end{equation}
We work with real solutions; the complex case follows by taking real and imaginary parts.

\subsection{Acknowledgments}
This work was partially supported by NSFC Grant (No.~12571260, 12422118), and the NSF of Hunan Province (No.~2026JJ20012).

{\bf AI disclosure.} During the preparation of this manuscript, the authors used ChatGPT to discuss possible proof strategies, assist in drafting and checking parts of the mathematical arguments, organize references, and improve the clarity and readability of the exposition. All mathematical statements and proofs in the final manuscript were written, reviewed, and verified by the authors, who take full responsibility for the content of this manuscript.

\section{Reduction to an observation estimate at one time}\label{sec-ob-one-time}

\subsection{Observation at one time and measurable-time telescoping}\label{subsec-time-telescoping}

The bounded-domain argument reduces to the following estimate. The energy
ratio $\Theta_{s,t}$ and its logarithmic parameter $\Lambda_{s,t}$ are defined
in the introduction.

\begin{theorem}[Observation estimate at one time]\label{thm-ob-one-time}
For every $m>0$ there exists $C_m=C(\Omega,T,M,m)$ such that, for
$0\le s<t\le T$ and every measurable $E\subset\Omega$ with $|E|\ge m$,
\begin{equation}\label{equ-111-1}
 \norm{u(t)}_2
 \le C_m\exp\!\left(C_m\Lambda_{s,t}^{3/2}\right)
                         \norm{u(t)}_{L^2(E)},
\end{equation}
whenever $\norm{u(t)}_2>0$.
\end{theorem}

Theorem~\ref{thm-ob-one-time} is proved in Subsection~\ref{subsec-proof-ob-one-time}.
We first establish the measurable-time argument that converts it into
Theorem~\ref{thm-ob-bound}. The same argument will be used for the whole-space
problem.

\begin{proposition}[Subquadratic observability criterion]\label{prop-subquadratic}
Let $h\in C([0,T])$ be nonnegative and nonincreasing, and let
$\mathcal A\subset(0,T)$ have positive measure. Suppose $g\in L^1(\mathcal A)$
is nonnegative and, for some $C_0>0$ and $1<\beta<2$,
\begin{equation}\label{equ-111-2}
 h(t)\le C_0\exp\!\left(C_0\left(1+(t-s)^{-1}
             +\log\frac{h(s)}{h(t)}\right)^\beta\right)g(t)
\end{equation}
for $0\le s<t\le T$, almost every $t\in\mathcal A$, whenever $h(t)>0$.
Then
\begin{equation}\label{equ-111-3}
 h(T)\le C\int_{\mathcal A}g(t)\,\mathrm{d}t,
 \qquad C=C(C_0,\beta,T,|\mathcal A|).
\end{equation}
The constant is uniform when $|\mathcal A|$ has a prescribed positive lower
bound.
\end{proposition}
\begin{proof}
\noindent\emph{Step 1. We introduce a small initial-state remainder.}

We first fix $0<\eta<1$. If $h(t)>\eta h(s)$, then $h(t)>0$ and
$\log(h(s)/h(t))<\log(1/\eta)$. It follows from \eqref{equ-111-2} that $h(t)$
is bounded by the observation term below. If $h(t)\le\eta h(s)$, the
initial-state term alone suffices. Combining these two cases, we obtain
\begin{equation}\label{equ-111-4}
 h(t)\le\eta h(s)+C\exp\!\left(
 C\bigl(1+(t-s)^{-1}+\log(1/\eta)\bigr)^\beta\right)g(t).
\end{equation}

\noindent\emph{Step 2. We choose intervals by observation measure.}

We define the accumulated observation time
\[
 \Phi(t)=|\mathcal A\cap(0,t)|,\qquad \kappa=|\mathcal A|/2,
\]
and choose
\begin{equation}\label{equ-111-5}
 1<p<1/(\beta-1),\qquad p\beta<\sigma<\beta/(\beta-1).
\end{equation}
For integers $k\ge1$, we set
\[
 m_k=\kappa k^{1-p},\qquad
 \tau_k=\inf\{t\in[0,T]:\Phi(t)\ge m_k\},\qquad
 \mu_k=m_k-m_{k+1}\approx k^{-p}.
\]
Here the comparison constants depend only on $\kappa,p$.
By construction, the function $\Phi$ is continuous and $1$-Lipschitz, with
$\Phi(\tau_k)=m_k$, and the sequence $\tau_k$ decreases. We take the latter half
of the observation measure in $(\tau_{k+1},\tau_k)$:
\[
 \mathcal A_k=\{t\in\mathcal A\cap(\tau_{k+1},\tau_k):
                         \Phi(t)-m_{k+1}\ge\mu_k/2\}.
\]
Since $\Phi'=\mathbf1_{\mathcal A}$ almost everywhere and $\Phi$ is
$1$-Lipschitz, we have
\begin{equation}\label{equ-111-6}
 |\mathcal A_k|=\mu_k/2,\qquad t-\tau_{k+1}\ge\mu_k/2
                                      \quad(t\in\mathcal A_k).
\end{equation}
By applying \eqref{equ-111-4} with $s=\tau_{k+1}$ and using
$h(\tau_k)\le h(t)$, we obtain, after averaging over $\mathcal A_k$,
\begin{equation}\label{equ-111-7}
 h(\tau_k)\le \eta h(\tau_{k+1}) +\frac C{\mu_k}
 \exp\!\left(C\bigl(1+\mu_k^{-1}+\log(1/\eta)\bigr)^\beta\right)
                         \int_{\mathcal A_k}g(t)\,\mathrm{d}t.
\end{equation}

\noindent\emph{Step 3. We choose summable weights.}

We set $w_k=\mathrm{e}^{-k^\sigma}$ and $\eta_k=w_{k+1}/w_k$. Then
$\log(1/\eta_k)\le Ck^{\sigma-1}$. Multiplying \eqref{equ-111-7}, with
$\eta=\eta_k$, by $w_k$, we find that the logarithm of the observation
coefficient is bounded above by
\[
 -k^\sigma+C\bigl(k^{p\beta}+k^{\beta(\sigma-1)}\bigr)+O(\log k).
\]
Both positive powers are strictly below $\sigma$ by
\eqref{equ-111-5}. Hence, for $k\ge k_0$, with $k_0$ depending only
on $C_0,\beta,T,|\mathcal A|$,
\[
 w_kh(\tau_k)-w_{k+1}h(\tau_{k+1})
 \le\int_{\mathcal A_k}g(t)\,\mathrm{d}t.
\]
By summing from $k_0$ to $K$ and letting $K\to\infty$, we obtain
\[
 w_{k_0}h(\tau_{k_0})\le\int_{\mathcal A}g(t)\,\mathrm{d}t,
\]
because $w_{K+1}h(\tau_{K+1})\le w_{K+1}h(0)\to0$ and the sets
$\mathcal A_k$ are disjoint. Finally, $h(T)\le h(\tau_{k_0})$, which proves
\eqref{equ-111-3}. Moreover, the construction is uniform for
$|\mathcal A|\ge a_0>0$, since $\kappa\in[a_0/2,T/2]$.
\end{proof}

\subsection{Bounded-domain observability}
\label{subsec-proof-ob-bound}

\begin{proof}[Proof of Theorem~\ref{thm-ob-bound} from Theorem~\ref{thm-ob-one-time}]
If $u(T)=0$, the conclusion is immediate. Otherwise, we put
\[
 m=\frac d{2T},\qquad
 \mathcal A=\{t\in(0,T):|D_t|\ge m\}.
\]
By Fubini's theorem, we have
\begin{equation}\label{equ-111-8}
 d\le\int_0^T|D_t|\,\mathrm{d}t
       \le mT+|\Omega|\,|\mathcal A|,
 \qquad |\mathcal A|\ge\frac d{2|\Omega|}.
\end{equation}
We set $I(t)=\int_{D_t}|u(x,t)|\,\mathrm{d}x$ and define
\begin{equation}\label{equ-111-9}
 \begin{aligned}
 h(t)&=\mathrm{e}^{-Mt}\norm{u(t)}_2,\qquad
 g(t)=\mathrm{e}^{-Mt}I(t),\\
 \frac{h(s)}{h(t)}&=\Theta_{s,t}^{1/2},\qquad
 \Lambda_{s,t}=1+(t-s)^{-1}+2\log\frac{h(s)}{h(t)}.
 \end{aligned}
\end{equation}
The energy regularity and energy inequality imply that $h$ is continuous and nonincreasing, with
$h(t)\ge h(T)>0$. Moreover, by Cauchy--Schwarz, we have
$0\le g(t)\le|\Omega|^{1/2}h(t)$, so $g\in L^1(\mathcal A)$.
For almost every $t\in\mathcal A$, we define
\begin{equation}\label{equ-112-1}
 F_t=\left\{x\in D_t:|u(x,t)|\le\frac{2I(t)}{|D_t|}\right\}.
\end{equation}
By applying Chebyshev's inequality to \eqref{equ-112-1}, we obtain
\begin{equation}\label{equ-112-2}
 \begin{aligned}
 |F_t|&\ge\frac{|D_t|}{2}\ge\frac m2,\\
 \norm{u(t)}_{L^2(F_t)}^2
 &\le\frac{2I(t)}{|D_t|}\int_{F_t}|u(x,t)|\,\mathrm{d}x
 \le\frac2m I(t)^2.
 \end{aligned}
\end{equation}
We now apply \eqref{equ-111-1} with $E=F_t$ and measure threshold $m/2$.
This, along with \eqref{equ-112-2}, gives
\begin{equation}\label{equ-112-3}
 \norm{u(t)}_2\le C\exp(C\Lambda_{s,t}^{3/2})I(t).
\end{equation}
By multiplying \eqref{equ-112-3} by $\mathrm{e}^{-Mt}$ and using
\eqref{equ-111-9}, we verify \eqref{equ-111-2}
with $\beta=3/2$, after adjusting the constant.
By \eqref{equ-111-3}, we therefore have
\begin{equation*}
 \norm{u(T)}_2
 \le C\mathrm{e}^{MT}\int_{\mathcal A}\mathrm{e}^{-Mt}I(t)\,\mathrm{d}t
 \le C\int_D|u(x,t)|\,\mathrm{d}x\,\mathrm{d}t.
\end{equation*}
The constant depends only on $\Omega,T,M,d$.
\end{proof}

It remains to establish the observation estimate at one time
\eqref{equ-111-1}. The following sections derive the small-ball mass bounds
and construct a cube on which the local observation estimate holds;
Subsection~\ref{subsec-proof-ob-one-time} then completes the proof.

\section{Local mass bounds and a cube criterion}\label{sec-local-estimates}
For this section and Section~\ref{sec-frequency-reduction}, fix \(0\le s<t\le T\)
with \(\norm{u(t)}_2>0\), and abbreviate
\(\Lambda=\Lambda_{s,t}\), \(\delta=t-s\).
Fix \(K\Subset\Omega\) and \(0<\rho<1\) such that
\(B_{8\rho}(x)\Subset\Omega\) for \(x\in K\).
All local constants are uniform over these centers.

\subsection{Gaussian frequency and small-ball mass}\label{subsec-local-frequency}
We use the localized Gaussian frequency method of \cite[Section 5]{PWZ}
and \cite[Section 2]{EFV}, and record the identities needed below.
For a smooth, spatially compactly supported \(F=F(y,\tau)\) and
\(\alpha>0\), put
\begin{equation}\label{equ-112-5}
 \begin{aligned}
 H_\alpha(\tau)&=\int_{\R^n}F^2G_{\alpha+\tau}\,\mathrm{d}y,\\
 D_\alpha(\tau)&=\int_{\R^n}|\nabla F|^2G_{\alpha+\tau}\,\mathrm{d}y,\\
 N_\alpha(\tau)&=\frac{2(\alpha+\tau)D_\alpha(\tau)}{H_\alpha(\tau)}.
 \end{aligned}
\end{equation}
Writing \(P=\partial_\tau+\Delta\), \cite[Lemma 2]{EFV} gives
\begin{align}
 H_\alpha'&=2D_\alpha+
       2\int_{\R^n}F(PF)G_{\alpha+\tau}\,\mathrm{d}y,
                                      \label{equ-112-6}\\
 N_\alpha'&\ge-\frac{\alpha+\tau}{H_\alpha}
       \int_{\R^n}|PF|^2G_{\alpha+\tau}\,\mathrm{d}y.
                                      \label{equ-112-7}
\end{align}
We also use the Gaussian moment inequality \cite[Lemma 3]{EFV}: for every \(\sigma>0\),
\begin{equation}\label{equ-112-8}
 \int_{\R^n}\frac{|y|^2}{8\sigma}F^2G_\sigma\,\mathrm{d}y
 \le2\sigma\int_{\R^n}|\nabla F|^2G_\sigma\,\mathrm{d}y
       +\frac n2\int_{\R^n}F^2G_\sigma\,\mathrm{d}y.
\end{equation}

\begin{lemma}\label{lem-ball-mass}
For all \(x\in K\),
\begin{equation}\label{equ-112-9}
 \int_{B_\rho(x)}u(y,t)^2\,\mathrm{d}y\ge \mathrm{e}^{-C\Lambda}\norm{u(t)}_2^2,\qquad
 \sup_{\tau\in[s,t]}\norm{u(\tau)}_2^2\le\Theta_{s,t}\norm{u(t)}_2^2.
\end{equation}
\end{lemma}
\begin{proof}
By a finite covering of $K$ by balls of radius $\rho/4$ and the interpolation
inequality \cite[Lemma 5]{PWZ}, applied on $[s,t]$ to the corresponding
balls of radius $\rho/2$, we obtain uniform $C>0$ and $0<\theta<1$ such that
\[
 \norm{u(t)}_2\le
 \mathrm{e}^{C(1+\delta^{-1})}
 \bigl(\mathrm{e}^{M\delta}\norm{u(s)}_2\bigr)^{1-\theta}
 \norm{u(t)}_{L^2(B_\rho(x))}^\theta,\qquad x\in K.
\]
By rearranging this inequality and using \eqref{equ-1031-7}, we obtain the
first estimate in \eqref{equ-112-9}. The second follows from the global
energy inequality.
\end{proof}

\begin{proposition}\label{prop-local-frequency}
For a standard cutoff $\chi_x$, equal to one on $B_{3\rho}(x)$ and supported
in $B_{4\rho}(x)$, there exist uniform $C,c>0$ such that
\begin{align}
 \frac{2r^2\int|\nabla(\chi_xu)(y,t-r^2)|^2G_{r^2}(y-x)\,\mathrm{d}y}
      {\int(\chi_xu)(y,t-r^2)^2G_{r^2}(y-x)\,\mathrm{d}y}
       &\le C\Lambda,\quad r\le c/\sqrt \Lambda,\label{equ-113-1}\\
 \frac{2\alpha\int|\nabla(\chi_xu)(y,t)|^2G_\alpha(y-x)\,\mathrm{d}y}
      {\int(\chi_xu)(y,t)^2G_\alpha(y-x)\,\mathrm{d}y}
       &\le C\Lambda,
       \quad\alpha\le c/\Lambda,\label{equ-113-2}\\
 \int_{B_{2r}(x)}u(y,t)^2\,\mathrm{d}y
       &\le \mathrm{e}^{C\Lambda}\int_{B_r(x)}u(y,t)^2\,\mathrm{d}y,
       \quad r\le c.\label{equ-113-3}
\end{align}
Moreover, for $g=(\partial_t-\Delta)(\chi_xu)$,
\begin{equation}\label{equ-113-4}
 \frac{\int g(y,t-r^2)^2G_{r^2}(y-x)\,\mathrm{d}y}
      {\int(\chi_xu)(y,t-r^2)^2G_{r^2}(y-x)\,\mathrm{d}y}\le C,
 \qquad r\le c/\sqrt \Lambda.
\end{equation}
\end{proposition}
\begin{proof}
For $q\in[t-\delta/2,t]$, the energy inequality and
\eqref{equ-1031-7} give $\Lambda_{s,q}\le2\Lambda$ and
$\norm{u(q)}_2\ge\mathrm{e}^{-MT}\norm{u(t)}_2$. By \eqref{equ-112-9}
at time $q$ and the energy and gradient estimates of
\cite[Lemma~1]{PWZ}, applied from time $s$, we obtain
\begin{align}
 \int_{B_\rho(x)}u(y,q)^2\,\mathrm{d}y
 &\ge\mathrm{e}^{-C\Lambda}\norm{u(t)}_2^2,\label{equ-113-5}\\
 \norm{u(q)}_{H^1(\Omega)}^2
 &\le C(1+\delta^{-1})\norm{u(s)}_2^2
 \le\mathrm{e}^{C\Lambda}\norm{u(t)}_2^2.\label{equ-113-6}
\end{align}

We put $F(y,\tau)=\chi_x(y)u(y,t-\tau)$, extend it by zero, and use
\eqref{equ-112-5} with Gaussians centered at $x$. Take $\beta=c_0/\Lambda\le\delta/4$,
$0\le\tau\le3\beta/4$, and $0<\alpha\le\beta/16$. With
$\sigma=\alpha+\tau$, we have $PF=a(y,t-\tau)F+R$, where
$R=2\nabla\chi_x\cdot\nabla u+(\Delta\chi_x)u$ at time $t-\tau$.
Since $R$ vanishes on $B_{3\rho}(x)$, \eqref{equ-113-6} bounds
$\int R^2G_\sigma(y-x)\,\mathrm{d}y$ by
$C\sigma^{-n/2}\exp(C\Lambda-9\rho^2/(4\sigma))\norm{u(t)}_2^2$.
Using \eqref{equ-113-5} on $B_\rho(x)$, we obtain
\begin{equation}\label{equ-113-7}
 \begin{aligned}
 H_\alpha(\tau)&\ge c\sigma^{-n/2}
 \exp\!\left(-C\Lambda-\frac{\rho^2}{4\sigma}\right)\norm{u(t)}_2^2,\\
 \frac{\int R^2G_\sigma(y-x)\,\mathrm{d}y}{H_\alpha(\tau)}
 &\le C\exp\!\left(C\Lambda-\frac{2\rho^2}{\sigma}\right)\le C
 \end{aligned}
\end{equation}
for $c_0$ sufficiently small. Thus $\int |PF|^2G_\sigma/H_\alpha\le C$,
and \eqref{equ-112-6}--\eqref{equ-112-7} give
\[
 N_\alpha'(\tau)\ge-C\sigma,\qquad
 \left|(\log H_\alpha)'(\tau)-\frac{N_\alpha(\tau)}{\sigma}\right|\le C.
\]
For $\tau\le\beta/4$ and $v\in[\beta/2,3\beta/4]$, the first inequality
gives $N_\alpha(\tau)\le N_\alpha(v)+C\beta^2$. Averaging in $v$ and
using $v+\alpha\le\beta$ in the second inequality, as in
\cite[(2.13)--(2.16)]{EFV}, yields
\begin{equation}\label{equ-113-8}
 \begin{aligned}
 N_\alpha(\tau)
 &\le C\beta^2+\frac4\beta\int_{\beta/2}^{3\beta/4}N_\alpha(v)\,\mathrm{d}v\\
 &\le C\left(1+\log\frac{H_\alpha(3\beta/4)}{H_\alpha(\beta/2)}\right)
 \le C(\Lambda+\beta^{-1})\le C\Lambda.
 \end{aligned}
\end{equation}
For the bound on the ratio, we used \eqref{equ-113-6} for the numerator
and \eqref{equ-113-7} for the denominator, with $\sigma$ comparable to
$\beta$ at both endpoints. Taking $\tau=0$ gives \eqref{equ-113-2}.
At fixed $\tau=r^2\in(0,\beta/4]$, letting $\alpha\downarrow0$ in
\eqref{equ-113-8} and in the bound on $\int |PF|^2G_\sigma/H_\alpha$
gives \eqref{equ-113-1} and \eqref{equ-113-4}, since $PF=-g$ at time
$t-r^2$.

To prove \eqref{equ-113-3}, we take $\alpha=r^2/(A\Lambda)$ with $A$
sufficiently large and $r\le c$ so that $\alpha\le\beta/16$ and
$B_{2r}(x)\subset B_{3\rho}(x)$. The moment inequality
\eqref{equ-112-8} and \eqref{equ-113-2} give
\[
 \int_{B_r(x)^c}F(y,0)^2G_\alpha(y-x)\,\mathrm{d}y
 \le\frac{8\alpha}{r^2}\bigl(N_\alpha(0)+n/2\bigr)H_\alpha(0)
 \le\tfrac12 H_\alpha(0).
\]
As in the proof of \cite[Lemma~4]{EFV}, comparison of the Gaussian on
$B_r(x)$ and $B_{2r}(x)$ then gives
\[
 \int_{B_{2r}(x)}u(y,t)^2\,\mathrm{d}y
 \le 2\exp(r^2/\alpha)\int_{B_r(x)}u(y,t)^2\,\mathrm{d}y.
\]
The frequency identities extend to energy solutions by approximation
as in \cite[Lemma~2]{PWZ}.
\end{proof}

\begin{lemma}\label{lem-small-ball-mass}
For allowed centers \(x\), \(t-\delta/2\le\tau\le t\), and \(0<\ell\le r_0\),
where \(r_0>0\) is fixed,
\begin{align}
 \int_{B_\ell(x)}u(y,\tau)^2\,\mathrm{d}y
 &\ge \mathrm{e}^{-C\Lambda(1+|\log\ell|)}\norm{u(t)}_2^2,\label{equ-113-9}\\
 \norm{u(\tau)}_2^2+\norm{\nabla u(\tau)}_2^2
 &\le \mathrm{e}^{C\Lambda}\norm{u(t)}_2^2.\label{equ-113-10}
\end{align}
\end{lemma}
\begin{proof}
For $\tau\in[t-\delta/2,t]$, we have $\Lambda_{s,\tau}\le2\Lambda$
and $\norm{u(\tau)}_2\ge\mathrm{e}^{-MT}\norm{u(t)}_2$ by the energy
inequality. We fix $0<r_0\le\min\{\rho,c\}$, with $c$ as in
\eqref{equ-113-3}. Let $m_\tau(r)$ be the squared $L^2$ mass on $B_r(x)$.
Applying \eqref{equ-112-9} at time $\tau$ with observation radius $r_0$
as in its proof using \cite[Lemma~5]{PWZ}, we have
$m_\tau(r_0)\ge\mathrm{e}^{-C\Lambda}\norm{u(t)}_2^2$.

For $r_j=2^{-j}r_0$, \eqref{equ-113-3} at time $\tau$ gives
$m_\tau(r_{j-1})\le\mathrm{e}^{C\Lambda}m_\tau(r_j)$ for $j\ge1$.
Taking $k=\lceil\log_2(r_0/\ell)\rceil$ and iterating, we obtain
\[
 m_\tau(\ell)\ge m_\tau(r_k)\ge\mathrm{e}^{-C\Lambda k}m_\tau(r_0)
 \ge\mathrm{e}^{-C\Lambda(k+1)}\norm{u(t)}_2^2.
\]
Since $k+1\le C(1+|\log\ell|)$, this proves \eqref{equ-113-9}.
Estimate \eqref{equ-113-10} follows directly from \eqref{equ-113-6}.
\end{proof}

\subsection{Common cutoffs and a cube criterion}\label{subsec-common-cutoff}
Choose \(\chi\in C_c^\infty(\Omega)\) equal to one on \(B_{2\rho}(x)\)
for every allowed center \(x\). Extend \(f=\chi u\) by zero to \(\R^n\),
and put \(g=(\partial_t-\Delta)f\).

\begin{lemma}\label{lem-common-cutoff}
For all allowed centers and \(0<r\le c/\sqrt\Lambda\),
\begin{equation}\label{equ-113-11}
 N_x(r):=\frac{2r^2\int|\nabla f(y,t-r^2)|^2G_{r^2}(y-x)\,\mathrm{d}y}
                   {\int f(y,t-r^2)^2G_{r^2}(y-x)\,\mathrm{d}y}
 \le C\Lambda,\qquad
 \frac{\int g(y,t-r^2)^2G_{r^2}(y-x)\,\mathrm{d}y}
      {\int f(y,t-r^2)^2G_{r^2}(y-x)\,\mathrm{d}y}\le C.
\end{equation}
\end{lemma}
\begin{proof}
Fix $q\in[t-\delta/2,t]$ and $0<\sigma\le c/\Lambda$. At time $q$,
let $v=\chi_xu$, and write $H(w)=\int w^2G_\sigma(y-x)\,\mathrm{d}y$
and $D(w)=\int|\nabla w|^2G_\sigma(y-x)\,\mathrm{d}y$. Both cutoffs
equal one on $B_{2\rho}(x)$, so \eqref{equ-113-5} on $B_\rho(x)$ gives
\[
 \min\{H(f),H(v)\}
 \ge c\sigma^{-n/2}\exp\!\left(-C\Lambda-\frac{\rho^2}{4\sigma}\right)
       \norm{u(t)}_2^2.
\]
Write $g=-af-R_\chi$, where
$R_\chi=2\nabla\chi\cdot\nabla u+(\Delta\chi)u$.
The differences of the integrands and the term $R_\chi$ vanish on
$B_{2\rho}(x)$. On its complement,
$G_\sigma\le C\sigma^{-n/2}\exp(-\rho^2/\sigma)$, so by
\eqref{equ-113-10} and the lower bound above,
\[
 \frac{|H(f)-H(v)|+|D(f)-D(v)|+\int R_\chi^2G_\sigma(y-x)\,\mathrm{d}y}
      {H(f)}
 \le C\exp\!\left(C\Lambda-\frac{3\rho^2}{4\sigma}\right)=:\varepsilon.
\]
After decreasing $c$, we have $\varepsilon\le1$. By
\eqref{equ-113-2} at time $q$ and $\Lambda_{s,q}\le2\Lambda$, we obtain
\[
 \begin{aligned}
 \frac{2\sigma D(f)}{H(f)}
 &\le (1+\varepsilon)\frac{2\sigma D(v)}{H(v)}+2\sigma\varepsilon
 \le C\Lambda,\\
 \frac{\int g^2G_\sigma(y-x)\,\mathrm{d}y}{H(f)}
 &\le 2M^2+2\varepsilon\le C.
 \end{aligned}
\]
Taking $q=t-r^2$ and $\sigma=r^2$ proves \eqref{equ-113-11}.
\end{proof}

\begin{proposition}\label{prop-cube-criterion}
Fix \(m>0\). Suppose that for every measurable \(E\subset\Omega\) with
\(|E|\ge m\) there is a cube \(Q=Q_\ell(x_Q)\), with allowed center and
\(0<\ell\le r_0\), such that
\begin{equation}\label{equ-114-1}
 \begin{aligned}
 |E\cap Q|&\ge c_m|Q|,\\
 \norm{u(t)}_{L^2(Q)}&\le C_m\norm{u(t)}_{L^2(E\cap Q)},\\
 1+|\log\ell|&\le C_m\sqrt\Lambda .
 \end{aligned}
\end{equation}
Then \eqref{equ-111-1} holds.
\end{proposition}
\begin{proof}
Since \(B_\ell(x_Q)\subset Q\), we obtain from \eqref{equ-113-9} at time \(t\) that
\[
 \norm{u(t)}_{L^2(Q)}^2
 \ge \mathrm{e}^{-C\Lambda(1+|\log\ell|)}\norm{u(t)}_2^2
 \ge \mathrm{e}^{-C_m\Lambda^{3/2}}\norm{u(t)}_2^2.
\]
Combining this with the second inequality in \eqref{equ-114-1}
and taking square roots, we obtain the desired estimate.
\end{proof}
To construct \(Q\), we first retain a positive observation density on a
cube whose past frequency has reached a fixed threshold at logarithmic
depth \(C_m\sqrt\Lambda\). The resulting frequency bound gives fixed terminal
cylinder doubling on a smaller scale. Propagation of smallness under this
fixed doubling bound gives the second inequality in \eqref{equ-114-1}.

\section{Frequency reduction and observation at one time}\label{sec-frequency-reduction}
\subsection{Gaussian variance and spatial covariance}

Let
\[
 W_{x,r}(y)=f(x+ry,t-r^2),\qquad \D=-2\Delta_y+y\cdot\nabla_y.
\]
On $L^2(\mathrm{d}\gamma)$ we use the nonnegative self-adjoint realization of $\D$.
Its degree-$k$ eigenspace is
\[
 \mathcal H_k=\{\mathrm{e}^{-\Delta}p:p\text{ is homogeneous of degree }k\},
 \qquad k=0,1,2,\ldots.
\]
Here $\mathrm{e}^{-\Delta}p=\sum_{j=0}^{\lfloor k/2\rfloor}(-\Delta)^jp/j!$ is a
finite polynomial. Let $\Pi_k$ be the orthogonal projection onto
$\mathcal H_k$. Then
\[
 \Dom(\D)=\left\{W\in L^2(\mathrm{d}\gamma):
                \sum_{k\ge0}k^2\norm{\Pi_kW}_\gamma^2<\infty\right\},
 \qquad \D W=\sum_{k\ge0}k\Pi_kW.
\]
This is the Gaussian caloric-degree decomposition used in \cite{HJ}.
For $W=W_{x,r}\in\Dom(\D)\setminus\{0\}$, define
\begin{equation}\label{equ-114-2}
 N_x(r)=\frac{\langle\D W,W\rangle_\gamma}{\norm W_\gamma^2},\qquad
 V_x(r)=\frac{\norm{(\D-N_x(r))W}_\gamma^2}{\norm W_\gamma^2}.
\end{equation}
The first expression equals the frequency in \eqref{equ-113-11}, which
is defined by its gradient formula at every scale. The variance is needed
only for almost every $r$, where $W_{x,r}\in\Dom(\D)$.
The numbers $\norm{\Pi_kW}_\gamma^2/\norm W_\gamma^2$ sum to one;
$N_x(r)$ and $V_x(r)$ are their mean degree and degree variance.

\begin{lemma}\label{lem-frequency-variance}
There exists a uniform $C_0$ such that
\begin{equation}\label{equ-114-3}
 r\frac{\mathrm{d}N_x}{\mathrm{d}r}\ge V_x-C_0r^4.
\end{equation}
Writing
$F_x(r)=N_x(r)+C_0r^4/4$, one has $rF_x'(r)\ge V_x(r)$.
Moreover, if the segment joining $x,y$ stays in the allowed interior region,
then, for $0<\tau_1<\tau_2\le c^2/\Lambda$,
\begin{equation}\label{equ-114-4}
 F_x(\sqrt{\tau_1})\le
       \exp\left(\frac{|x-y|^2}{\tau_2-\tau_1}\right)F_y(\sqrt{\tau_2}).
\end{equation}
\end{lemma}
\begin{proof}
We fix the center $x$ and write $W=W_{x,r}$, $N=N_x(r)$, and $V=V_x(r)$.
We set $g_r(y)=g(x+ry,t-r^2)$.
By changing variables and applying Lemma~\ref{lem-common-cutoff}, we obtain the uniform bound
\[
 \frac{\norm{g_r}_\gamma^2}{\norm W_\gamma^2}
 =\frac{\int_{\R^n}|g(z,t-r^2)|^2G_{r^2}(z-x)\,\mathrm{d}z}
        {\int_{\R^n}|f(z,t-r^2)|^2G_{r^2}(z-x)\,\mathrm{d}z}
 \le C_*.
\]
With $x,t$ fixed, the chain rule and $\partial_t f=\Delta f+g$ give
\[
 \begin{aligned}
 r\partial_rW(y)
 &=r y\cdot(\nabla f)(x+ry,t-r^2)
       -2r^2(\partial_t f)(x+ry,t-r^2)\\
 &=y\cdot\nabla_yW-2\Delta_yW-2r^2g_r
  =\D W-2r^2g_r.
 \end{aligned}
\]
We set $Z=(\D-N)W$. By the definition of $N$, we have
$\langle Z,W\rangle_\gamma=0$, and consequently
$\langle Z,\D W\rangle_\gamma=\norm Z_\gamma^2$.
Since $\D$ is self-adjoint, we differentiate the numerator and denominator
of the Rayleigh quotient to obtain, at almost every admissible scale,
\[
 \begin{aligned}
 r\frac{\mathrm{d}N}{\mathrm{d}r}
 &=\frac{2\langle\D W,r\partial_rW\rangle_\gamma}{\norm W_\gamma^2}
   -\frac{2N\langle W,r\partial_rW\rangle_\gamma}{\norm W_\gamma^2}\\
 &=\frac{2\langle Z,\D W-2r^2g_r\rangle_\gamma}{\norm W_\gamma^2} =2V-4r^2\frac{\langle Z,g_r\rangle_\gamma}{\norm W_\gamma^2}.
 \end{aligned}
\]
By the source bound above, Cauchy--Schwarz, and Young's inequality, we obtain
\[
 4r^2\frac{|\langle Z,g_r\rangle_\gamma|}{\norm W_\gamma^2}
 \le4\sqrt{C_*}\,r^2\sqrt V\le V+4C_*r^4.
\]
It follows that $rN'(r)\ge V-4C_*r^4$.
By taking $C_0=4C_*$, we obtain \eqref{equ-114-3} and
\begin{equation}\label{equ-114-5}
 rF_x'(r)=rN_x'(r)+C_0r^4\ge V_x(r).
\end{equation}
Differentiation with respect to the center and
$2\norm{\nabla_yW}_\gamma^2=N\norm W_\gamma^2$ give
\begin{equation}\label{equ-114-6}
 \begin{aligned}
 \partial_{x_i}N_x(r)&=\frac2r
       \frac{\langle(\D-N)W,\partial_{y_i}W\rangle_\gamma}{\norm W_\gamma^2},\\
 |\nabla_xN_x(r)|^2&\le\frac{2N_x(r)V_x(r)}{r^2}.
 \end{aligned}
\end{equation}
We now set $\mathcal F(x,\tau)=F_x(\sqrt\tau)$.
By \eqref{equ-114-5} and
\eqref{equ-114-6}, using $N_x\le F_x$, we have
\begin{equation}\label{equ-114-7}
 \partial_\tau\mathcal F\ge\frac{V_x(\sqrt\tau)}{2\tau},\qquad
 |\nabla_x\mathcal F|^2\le4\mathcal F\,\partial_\tau\mathcal F.
\end{equation}
For $\varepsilon>0$, we set $\Phi_\varepsilon=\log(\mathcal F+\varepsilon)$.
Since $\mathcal F\ge0$, \eqref{equ-114-7} gives
$\partial_\tau(\Phi_\varepsilon) \ge|\nabla\Phi_\varepsilon|^2/4$.
Along the straight path
$z(\tau)=x+(\tau-\tau_1)(y-x)/(\tau_2-\tau_1)$, we obtain
\[
 \frac{\mathrm{d}}{\mathrm{d}\tau}\Phi_\varepsilon(z(\tau),\tau)
 \ge\tfrac14|\nabla\Phi_\varepsilon|^2
                  +z'\cdot\nabla\Phi_\varepsilon\ge-|z'|^2.
\]
By integrating from $\tau_1$ to $\tau_2$, we obtain
\[
 \Phi_\varepsilon(y,\tau_2)-\Phi_\varepsilon(x,\tau_1)
       \ge-\frac{|x-y|^2}{\tau_2-\tau_1}.
\]
Exponentiating, we obtain
\[
 \mathcal F(x,\tau_1)+\varepsilon
 \le\exp\!\left(\frac{|x-y|^2}{\tau_2-\tau_1}\right)
                    (\mathcal F(y,\tau_2)+\varepsilon).
\]
By letting $\varepsilon\downarrow0$, we obtain \eqref{equ-114-4}.
\end{proof}

\begin{lemma}\label{lem-hermite}
Let $v\in L^2(\R^n)$ satisfy $(-\Delta+|x|^2)v=\lambda^2v$ in $\R^n$,
with $\lambda>0$. For $R\ge1$ and $|z|\le\lambda+1$,
\begin{equation}\label{equ-114-8}
 \norm v_{L^2(B_R(z))}^2
       \le C_n(R/\lambda)^{1/2}\norm v_{L^2(\R^n)}^2.
\end{equation}
\end{lemma}
\begin{proof}
We may assume $v\ne0$.
For $|z|\le\lambda$, we set
$\eta=\max\{\lambda^{-4/3},1-|z|/\lambda\}$.
By \cite[Theorem~1 and (1.16) in the arXiv version]{WZ}, for $0<R\le\lambda$ we have
\begin{equation}\label{equ-114-9}
 \frac{\norm v_{L^2(B_R(z))}}{\norm v_{L^2(\R^n)}}
 \le C_n
 \begin{cases}
  (\lambda\eta^{1/2})^{(n-2)/2}R^{n/2},
       & R\lesssim(\lambda\eta^{1/2})^{-1},\\
  R^{1/2}\lambda^{-1/2}\eta^{-1/4},
       & (\lambda\eta^{1/2})^{-1}\ll R\ll\lambda\eta,\\
  (R/\lambda)^{1/4},&\lambda\eta\lesssim R\le\lambda.
 \end{cases}
\end{equation}
For $R\ge1$ and sufficiently large $\lambda$, the first range is excluded
since $(\lambda\eta^{1/2})^{-1}\le\lambda^{-1/3}$.
In the second range, it follows from $R\lesssim\lambda\eta$ that
\[
 R^{1/2}\lambda^{-1/2}\eta^{-1/4}
       \le C(R/\lambda)^{1/4};
\]
the third range gives the same bound directly.
The trivial $L^2$ bound gives the result for bounded $\lambda$ and $R\gtrsim\lambda$.
For $\lambda<|z|\le\lambda+1$, we take $z_0=\lambda z/|z|$.
Since $B_R(z)\subset B_{R+1}(z_0)$ and $R+1\le2R$, the preceding bound
still applies. Thus we obtain \eqref{equ-114-8}.
\end{proof}

Let $0\ne W\in L^2(\mathrm{d}\gamma)$ have finite second moment,
$\int_{\R^n}|y|^2|W(y)|^2\,\mathrm{d}\gamma(y)<\infty$.
Define the probability measure and its mean vector by
\begin{equation}\label{equ-114-10}
 \mathrm{d}\mu_W(y)=\frac{|W(y)|^2}{\norm W_\gamma^2}\,\mathrm{d}\gamma(y),
 \qquad \bar y_W=\int_{\R^n}y\,\mathrm{d}\mu_W(y).
\end{equation}
Its covariance matrix and covariance trace are
\begin{align}
 \Cov(\mu_W)&=\int_{\R^n}
       (y-\bar y_W)\otimes(y-\bar y_W)\,\mathrm{d}\mu_W(y),\label{equ-115-1}\\
 S(W)&=\operatorname{tr}\Cov(\mu_W)
       =\int_{\R^n}|y-\bar y_W|^2\,\mathrm{d}\mu_W(y).
                                                        \label{equ-115-2}
\end{align}
Here $p\otimes p$ denotes the matrix with entries $p_ip_j$.
When $W$ is fixed, we abbreviate $\bar y_W$ to $\bar y$.
The following estimate connects the degree distribution to this spatial
distribution.

\begin{lemma}\label{lem-degree-covariance}
For $0\ne W\in\Dom(\D)$, let
\[
 N=\frac{\langle\D W,W\rangle_\gamma}{\norm W_\gamma^2},\qquad
 V=\frac{\norm{(\D-N)W}_\gamma^2}{\norm W_\gamma^2},
 \qquad S=S(W).
\]
Then
\begin{equation}\label{equ-115-3}
 \frac{N}{(1+V)^2}\le C_n(1+S).
\end{equation}
\end{lemma}
\begin{proof}
We normalize $\norm W_\gamma=1$. The unitary map $\mathcal U$ given by
\[
 w(x)=(\mathcal UW)(x)=\pi^{-n/4}\mathrm{e}^{-|x|^2/2}W(2x)
\]
conjugates $\D$ to $(-\Delta_x+|x|^2-n)/2$. Thus degree $k$ corresponds
to oscillator eigenvalue $2k+n$. The probability density
$|w(x)|^2\,\mathrm{d}x$ has mean and centered second moment
\[
 \bar x=\int_{\R^n}x|w(x)|^2\,\mathrm{d}x=\bar y/2,\qquad
 \int_{\R^n}|x-\bar x|^2|w(x)|^2\,\mathrm{d}x=S/4.
\]
In particular,
$|\bar x|^2\le\int|x|^2|w|^2\le
\langle(-\Delta+|x|^2)w,w\rangle=2N+n$.

The conclusion is immediate if $1+V\ge\sqrt N$ or $N$ is bounded. Otherwise, we set
\[
 J=8\sqrt{1+V}\le8N^{1/4},\qquad
 I=\{k\in\mathbb N_0:|k-N|\le J\},\qquad
 w_I=\sum_{k\in I}\mathcal U\Pi_kW.
\]
Here $\mathbb N_0=\{0,1,2,\ldots\}$. By Chebyshev's inequality for the degree distribution, we obtain
\[
 \begin{aligned}
 \norm{w-w_I}_2^2
 &=\sum_{k\notin I}\norm{\Pi_kW}_\gamma^2\\
 &\le\frac1{J^2}\sum_{k\ge0}(k-N)^2\norm{\Pi_kW}_\gamma^2
 =\frac V{J^2}\le\frac1{64}.
 \end{aligned}
\]
Applying Chebyshev's inequality in space with $R=2\sqrt{1+S}$, we also obtain
\[
 \int_{B_R(\bar x)^c}|w|^2
 \le\frac1{R^2}\int_{\R^n}|x-\bar x|^2|w(x)|^2\,\mathrm{d}x
 =\frac{S}{4R^2}\le\frac1{16}.
\]
The triangle inequality therefore yields
\begin{equation}\label{equ-115-4}
 \norm{w_I}_{L^2(B_R(\bar x))}
 \ge\norm w_{L^2(B_R(\bar x))}-\norm{w-w_I}_2
 \ge\sqrt{15/16}-1/8\ge1/2.
\end{equation}

For every degree in this band, $N/2\le k\le2N$ and
\[
 |\bar x|\le\sqrt{2N+n}\le\sqrt{2k+n}+CJ/\sqrt N
                    \le\sqrt{2k+n}+1
\]
when $N$ is sufficiently large. By Lemma~\ref{lem-hermite},
Cauchy--Schwarz over the $O(J+1)$ degrees, and their global orthogonality, we obtain
\[
 \tfrac14\le\norm{w_I}_{L^2(B_R(\bar x))}^2
       \le C_n(J+1)(R/\sqrt N)^{1/2}\norm{w_I}_2^2.
\]
Consequently, $R\ge c_n\sqrt N/(J+1)^2$, which yields \eqref{equ-115-3}.
\end{proof}

\begin{lemma}\label{lem-mean-covariance}
Let $Q$ be a cube of half-side $R$, let $\lambda\ge1$, and suppose
$N_x(r)\le\lambda$ for $x\in\overline Q$ and $r\le R/\sqrt\lambda$.
Let $\mu_{x,r}=\mu_{W_{x,r}}$ be the probability measure defined by
\eqref{equ-114-10}, with mean vector
$\bar y_{x,r}=\int_{\R^n}y\,\mathrm{d}\mu_{x,r}(y)$. Define
\begin{equation}\label{equ-115-5}
 S_x(r)=\operatorname{tr}\Cov(\mu_{x,r})
       =\int_{\R^n}|y-\bar y_{x,r}|^2\,\mathrm{d}\mu_{x,r}(y).
\end{equation}
Then
\begin{equation}\label{equ-115-6}
 \avg_Q S_x(r) \,\mathrm{d}x\le C_n.
\end{equation}
\end{lemma}
\begin{proof}
We fix the scale $r$ and write
\[
 H(x,r)=\norm{W_{x,r}}_\gamma^2
       =\int_{\R^n}f(z,t-r^2)^2G_{r^2}(z-x)\,\mathrm{d}z.
\]
By differentiating with respect to the physical center $x$, we obtain
\[
 \nabla_x\log H=\frac{\bar y_{x,r}}{2r},\qquad
 4r^2\nabla_x^2\log H=\Cov(\mu_{x,r})-2I_n,
\]
where $I_n$ is the $n\times n$ identity matrix. By taking the trace, we find
$S_x(r)=4r^2\Delta_x\log H+2n$.
By integration by parts in $z$, followed by Cauchy--Schwarz, we obtain
\[
 \begin{aligned}
 \nabla_x H&=2\int_{\R^n}f(z,t-r^2)\nabla_z f(z,t-r^2)
                         G_{r^2}(z-x)\,\mathrm{d}z,\\
 |\nabla_x\log H|
 &\le2\left(\frac{\int_{\R^n}|\nabla_z f(z,t-r^2)|^2
                         G_{r^2}(z-x)\,\mathrm{d}z}{H}\right)^{1/2}
 =\frac{\sqrt{2N_x(r)}}r.
 \end{aligned}
\]
Since $N_x(r)\le\lambda$ on $\partial Q$, the divergence theorem gives
\[
 \begin{aligned}
 \avg_QS_x(r)\,\mathrm{d}x
 &=2n+\frac{4r^2}{|Q|}\int_{\partial Q}
                 \nabla_x\log H\cdot\nu\,\mathrm{d}\sigma\\
 &\le2n+4n\sqrt2\,\frac{r\sqrt\lambda}{R}\le C_n.
 \end{aligned}
\]
Here $\nu$ is the outward unit normal and $\mathrm{d}\sigma$ is surface
measure. This proves \eqref{equ-115-6}.
\end{proof}

\subsection{Collective frequency reduction and covering}\label{subsec-frequency-covering}

The previous two lemmas have complementary roles: large degree with small
degree variance forces large spatial covariance, while the covariance
averaged over centers stays bounded. Thus high frequency can persist at
many centers only if a definite amount of variance is spent as the scale changes.

\begin{proposition}[Collective frequency halving]\label{prop-frequency-halving}
There is a constant $A_n\ge1$, depending only on $n$, with the following
property.
Let $Q$ have half-side $R$, and suppose the common-cutoff framework holds on
$3Q$, with
\[
 F_x(r)\le\lambda\quad(x\in3Q,\ 0<r\le r_{\max}),\qquad
 C_0r_{\max}^4/4\le1/10.
\]
For $0<\epsilon<1/4$ and $\epsilon\lambda\ge A_n$, set
\begin{equation}\label{equ-116-1}
 a=\min\{r_{\max}/2,R/\sqrt\lambda\},\qquad
 b=a\exp\!\left(-2\sqrt{A_n\lambda/\epsilon}\right).
\end{equation}
Then
\begin{equation}\label{equ-116-2}
 |\{x\in Q:F_x(b)>\lambda/2\}|\le\epsilon|Q|.
\end{equation}
\end{proposition}
\begin{proof}
Since $a\le R/\sqrt\lambda$ and $N_x(r)\le F_x(r)\le\lambda$,
Lemmas~\ref{lem-degree-covariance} and \ref{lem-mean-covariance} allow us
to fix $A_n\ge1$, depending only on $n$, so that for almost every $r\in(0,a]$,
\begin{equation}\label{equ-1111-1}
 \int_Q\frac{N_x(r)}{(1+V_x(r))^2}\,\mathrm{d}x\le\frac{A_n}{54}|Q|.
\end{equation}
The hypothesis $\epsilon\lambda\ge A_n$ gives
$v:=2\sqrt{\epsilon\lambda/A_n}\ge2$. On
\[
 E_r=\{x\in Q:N_x(r)\ge\lambda/3,\ V_x(r)<v\},
\]
the integrand in \eqref{equ-1111-1} is at least $\lambda/[3(1+v)^2]$.
Since $1+v\le3v/2$, we obtain
\begin{equation}\label{equ-1111-2}
 |E_r|\le\frac{A_n(1+v)^2}{18\lambda}|Q|
       \le\frac{A_n v^2}{8\lambda}|Q|=\frac\epsilon2|Q|.
\end{equation}
Suppose $B=\{x\in Q:F_x(b)>\lambda/2\}$ has measure greater than $\epsilon|Q|$.
For $x\in B$ and $b\le r\le a$, monotonicity and $F_x-N_x\le1/10$ give
$N_x(r)\ge\lambda/2-1/10\ge\lambda/3$, since $\lambda>4$.
Thus $V_x(r)\ge v$ on $B\setminus E_r$. By \eqref{equ-1111-2}, for almost
every $r\in[b,a]$,
\begin{equation}\label{equ-1111-3}
 \int_BV_x(r)\,\mathrm{d}x
 \ge v\bigl(|B|-|E_r|\bigr)\ge\frac{v|B|}{2}.
\end{equation}
Integrating $rF_x'(r)\ge V_x(r)$ and using $0\le F_x\le\lambda$ and
\eqref{equ-1111-3}, we obtain
\[
 \begin{aligned}
 \lambda|B|&\ge\int_B\bigl(F_x(a)-F_x(b)\bigr)\,\mathrm{d}x\ge\int_b^a\left(\int_BV_x(r)\,\mathrm{d}x\right)\frac{\mathrm{d}r}{r}\\
 &\ge\frac{v|B|}{2}\log(a/b)=2\lambda|B|.
 \end{aligned}
\]
The last equality follows from \eqref{equ-116-1}. This contradiction
proves \eqref{equ-116-2}.
\end{proof}

\begin{proposition}
\label{prop-frequency-covering}
Fix an interior cube $Q_0$ of half-side $R_0<1$, with $3Q_0$ in the allowed
center region. If $E\subset Q_0$ and $|E|\ge\mu|Q_0|>0$, there are a fixed
$\lambda_* =\lambda_*(n,\mu)$ and a subcube $Q_*$ of half-side $R_*$ such that
\begin{align}
 \exp(-C_\mu \Lambda^{1/2})&\le R_*\le c/\sqrt \Lambda,\label{equ-116-3}\\
 |E\cap Q_*|&\ge(\mu/2)|Q_*|,\label{equ-116-4}\\
 F_x(r)&\le\lambda_*\quad(x\in3Q_*,\ 0<r\le R_*).
                                                        \label{equ-116-5}
\end{align}
The constants may depend on the fixed interior geometry.
\end{proposition}
\begin{proof}
\noindent\emph{Step 1. We choose the parameters for the iteration.}

By Lemma~\ref{lem-common-cutoff}, the initial bound is
$F_x(r)\le\lambda_0$ for $x\in3Q_0$ and $0<r\le r_{\max,0}$, where
\[
 \lambda_0=C\Lambda,\qquad r_{\max,0}=c/\sqrt\Lambda,\qquad
 C_0r_{\max,0}^4/4\le1/10
\]
after decreasing $c$. We set $\vartheta=3/4$ and
$\lambda_j=\vartheta^j\lambda_0$, and let $A_n$ be as in
Proposition~\ref{prop-frequency-halving}.
We choose $\eta=\eta(\mu)>0$ and $\lambda_*=\lambda_*(n,\mu)$ so that
\[
 \frac{\eta}{1-\vartheta^{1/3}}
 +\frac{A_n}{(1-\vartheta)\lambda_*}<\mu/4,
 \qquad \eta+A_n/\lambda_*<1/4.
\]
Let $J$ be the first index with $\lambda_J\le\lambda_*$.
We first consider the case $J\ge1$ and define
\begin{equation}\label{equ-116-6}
 \epsilon_j=\eta(\lambda_j/\lambda_0)^{1/3}+A_n/\lambda_j,
 \qquad 0\le j<J.
\end{equation}
Since $\lambda_{J-1}>\lambda_*$ and
$\lambda_j^{-1}=\vartheta^{J-1-j}/\lambda_{J-1}$, we obtain
\[
 \begin{gathered}
 0<\epsilon_j<1/4,\qquad \epsilon_j\lambda_j\ge A_n,\\
 \sum_{j<J}\epsilon_j
 \le\frac{\eta}{1-\vartheta^{1/3}}
   +\frac{A_n}{(1-\vartheta)\lambda_*}<\mu/4.
 \end{gathered}
\]

\noindent\emph{Step 2. We construct a cube satisfying \eqref{equ-116-4} and \eqref{equ-116-5}.}

Starting with $Q_0$, we construct the retained cubes by induction.
Suppose that each retained cube $Q$ at level $j$ has half-side $R_j$ and satisfies
\[
 F_x(r)\le\lambda_j\qquad(x\in3Q,\ 0<r\le r_{\max,j}).
\]
We apply Proposition~\ref{prop-frequency-halving} to each $Q$ with $\epsilon_j$.
The scales $a_j,b_j$ in \eqref{equ-116-1} are common to the whole level.
We subdivide each $Q$ dyadically into cubes of half-side $R_{j+1}$ satisfying
\[
 c_1b_j\le R_{j+1}\le c_2b_j,\qquad c_1=c_2/2,\quad 0<c_2\le1/4.
\]
We retain a child $P\subset Q$ if it contains $y$ with
$F_y(b_j)\le\lambda_j/2$. For $x\in3P\subset3Q$, we have
$|x-y|\le C_nR_{j+1}$. By \eqref{equ-114-4}, we obtain
\[
 \begin{aligned}
 F_x(b_j/2)
 &\le\exp\!\left(\frac{|x-y|^2}{3b_j^2/4}\right)F_y(b_j)\\
 &\le\exp(C_nc_2^2)\frac{\lambda_j}{2}
 \le\frac34\lambda_j=\lambda_{j+1},
 \end{aligned}
\]
where $c_2$ is sufficiently small. By monotonicity, this gives the
required bound at level $j+1$, with $r_{\max,j+1}=b_j/2\ge R_{j+1}$.

Let $U_j$ be the union of the retained level-$j$ cubes, with $U_0=Q_0$.
Since each discarded child lies in the exceptional set in
\eqref{equ-116-2}, we have
\[
 |U_j\setminus U_{j+1}|
 \le\epsilon_j|U_j|\le\epsilon_j|Q_0|.
\]
It follows that
\[
 \begin{aligned}
 |E\cap U_J|
 &\ge |E|-\sum_{j<J}|U_j\setminus U_{j+1}|\\
 &\ge\frac{3\mu}{4}|Q_0|\ge\frac{3\mu}{4}|U_J|.
 \end{aligned}
\]
By averaging over the level-$J$ cubes, we can find a cube $Q_*$
satisfying \eqref{equ-116-4}. Since $\lambda_J\le\lambda_*$,
we also obtain \eqref{equ-116-5}.

\noindent\emph{Step 3. We prove \eqref{equ-116-3}.}

For $1\le j<J$,
\[
 r_{\max,j}/2=b_{j-1}/4\ge R_j,\qquad a_j=R_j/\sqrt{\lambda_j}.
\]
At $j=0$, fixed geometry gives $a_0\ge cR_0/\sqrt{\lambda_0}$.
By \eqref{equ-116-1} and $R_{j+1}\ge c_1b_j$, we obtain
\[
 \log(R_j/R_{j+1})
 \le C\sqrt{\lambda_j/\epsilon_j}+C(1+\log\lambda_j).
\]
By \eqref{equ-116-6} and the definition of $J$, we have
\[
 \sqrt{\lambda_j/\epsilon_j}
 \le\eta^{-1/2}\sqrt{\lambda_0}\,\vartheta^{j/3},
 \qquad J\le C\log(2+\Lambda).
\]
Combining these estimates and summing over $j$, we get
\[
 \begin{aligned}
 \log(R_0/R_*)
 &\le C\eta^{-1/2}\sqrt{\lambda_0}\sum_{j<J}\vartheta^{j/3}
       +C\log^2(2+\Lambda)\\
 &\le C_\mu\sqrt\Lambda.
 \end{aligned}
\]
This proves the lower bound in \eqref{equ-116-3}. The upper bound follows
from $R_*\le R_1\le c_2b_0\le r_{\max,0}=c/\sqrt\Lambda$.

Finally, if $J=0$, then $\lambda_0\le\lambda_*$ and $r_{\max,0}$ has a fixed
positive lower bound. We subdivide $Q_0$ at a fixed dyadic scale below
this bound and choose a child with observation density at least $\mu$
by averaging. The initial frequency bound applies on its triple, and
\eqref{equ-116-3} follows after adjusting $C_\mu$.
\end{proof}

\subsection{From past frequency to terminal doubling}\label{subsec-terminal-doubling}

The following proposition uses the frequency bound on \(3Q_*\) from
Proposition~\ref{prop-frequency-covering} to obtain doubling on cylinders
ending at \(t\). The doubling constant is independent of \(\Lambda\),
as required in the application of Proposition~\ref{prop-propagation}.

\begin{proposition}[Terminal cylinder doubling]
\label{prop-terminal-doubling}
Let $Q_*$ satisfy \eqref{equ-116-3} and
\eqref{equ-116-5}, with fixed $\lambda_*$.
There are constants $C_\mu,D_*>0$, independent of $\Lambda$, such that
\begin{equation}\label{equ-116-7}
 \avg_{\Q_{2r}(x,t)}u^2\le4^{D_*}\avg_{\Q_r(x,t)}u^2,
 \qquad x\in\tfrac32Q_*,\quad
 0<r\le\frac{R_*}{C_\mu(1+\Lambda)^2}.
\end{equation}
\end{proposition}
\begin{proof}
By the hypotheses, we have $N_x(r)\le\lambda_*$ for $x\in3Q_*$ and $r\le R_*$.
For a fixed constant $A_\mu$ to be chosen, we set
\begin{equation}\label{equ-116-8}
 h=R_*/[A_\mu(1+\Lambda)^2],
\end{equation}
so it suffices to prove \eqref{equ-116-7} for $r\le ch$, after which we take
$C_\mu=A_\mu/c$.
We write
\[
 H_\sigma(x,q)=\int_{\R^n}f(y,q)^2G_\sigma(y-x)\,\mathrm{d}y.
\]

For $x\in\frac32Q_*$, $0<\alpha\le h^2$, and $0\le\tau\le h^2$, we claim that
\begin{equation}\label{equ-116-9}
 N_\alpha^x(\tau):=
 \frac{2(\alpha+\tau)\int|\nabla f(y,t-\tau)|^2G_{\alpha+\tau}(y-x)\,\mathrm{d}y}
      {H_{\alpha+\tau}(x,t-\tau)}
       \le C(\lambda_*+1).
\end{equation}
We take $A_\mu$ large enough that $h^2\le\delta/2$ and
$2h^2\le\min\{c/\Lambda,r_0^2\}$, with $r_0$ from
Lemma~\ref{lem-small-ball-mass}. Writing $\sigma=\alpha+\tau$ and $q=t-\tau$,
we apply the source estimate in the proof of Lemma~\ref{lem-common-cutoff}
at $(q,\sigma)$ to obtain
\begin{equation}\label{equ-116-10}
 \frac{\int g(y,t-\tau)^2G_{\alpha+\tau}(y-x)\,\mathrm{d}y}
      {H_{\alpha+\tau}(x,t-\tau)}
 \le CM^2+C\mathrm{e}^{C\Lambda-c/(\alpha+\tau)}\le C,
\end{equation}
where $g=(\partial_t-\Delta)f$.
We first consider $\alpha\le\tau\le h^2$. By
$G_\alpha*G_\tau=G_{\alpha+\tau}$, we obtain
\begin{equation}\label{equ-116-11}
 N_\alpha^x(\tau)=\frac{\alpha+\tau}{\tau}
 \frac{\int N_z(\sqrt\tau)H_\tau(z,t-\tau)G_\alpha(z-x)\,\mathrm{d}z}
      {H_{\alpha+\tau}(x,t-\tau)}.
\end{equation}
The integral over $3Q_*$ in \eqref{equ-116-11} is at most
$\lambda_* H_\sigma(x,q)$ by \eqref{equ-116-5} and the Gaussian convolution
identity. For $z\notin3Q_*$, we have $|z-x|\ge R_*$ and hence
$G_\alpha(z-x)\le C\exp(-R_*^2/(8\alpha))G_{2\alpha}(z-x)$.
Using $N_z(\sqrt\tau)H_\tau(z,q)=2\tau\int|\nabla f(y,q)|^2
G_\tau(y-z)\,\mathrm{d}y$, convolution and \eqref{equ-113-10} give
\[
 \int_{(3Q_*)^c}N_z(\sqrt\tau)H_\tau(z,q)G_\alpha(z-x)\,\mathrm{d}z
 \le C\tau\sigma^{-n/2}\exp\!\left(C\Lambda-\frac{R_*^2}{8\alpha}\right)
       \norm{u(t)}_2^2.
\]
By \eqref{equ-113-9} on $B_{\sqrt\sigma}(x)$, we have
\[
 H_\sigma(x,q)\ge c\sigma^{-n/2}
 \exp(-C\Lambda(1+|\log\sqrt\sigma|))\norm{u(t)}_2^2.
\]
Their ratio is thus bounded by
\[
 C\tau\exp\!\left(C\Lambda(1+|\log\sqrt\sigma|)-\frac{R_*^2}{8\alpha}\right)
 \le1
\]
for sufficiently large $A_\mu$. Indeed, with $b=h^2/\alpha\ge1$,
\eqref{equ-116-3} and \eqref{equ-116-8} give
\[
 |\log\sqrt\sigma|\le|\log h|+\tfrac12\log b
 \le C_\mu\sqrt\Lambda+\log A_\mu+\tfrac12 b,\qquad
 R_*^2/\alpha=A_\mu^2(1+\Lambda)^4 b.
\]
Since $(\alpha+\tau)/\tau\le2$, \eqref{equ-116-11} gives
$N_\alpha^x(\tau)\le2(\lambda_*+1)$.

For $0\le\tau<\alpha$, \eqref{equ-116-10} and \eqref{equ-112-7} give
\[
 (N_\alpha^x)'(\tau)\ge-C(\alpha+\tau).
\]
Integrating up to $\tau=\alpha$ and using the preceding bound, we obtain
\[
 N_\alpha^x(\tau)
 \le N_\alpha^x(\alpha)+C\int_\tau^\alpha(\alpha+s)\,\mathrm{d}s
 \le C(\lambda_*+1).
\]
This proves \eqref{equ-116-9}.

To prove \eqref{equ-116-7}, we follow the passage from Gaussian mass
to cylinder mass in \cite[proof of Theorem~3, (2.19)--(2.22)]{EFV}.
By the energy identity \eqref{equ-112-6},
\eqref{equ-116-9}, and \eqref{equ-116-10}, we have
\[
 \frac{\mathrm{d}}{\mathrm{d}\tau}
 \log H_{\alpha+\tau}(x,t-\tau)
 \le\frac{C(\lambda_*+1)}{\alpha+\tau}+C.
\]
Integration from $0$ to $\tau$ gives
\begin{equation}\label{equ-117-1}
 H_{\alpha+\tau}(x,t-\tau)\le
 \mathrm{e}^{C\tau}(1+\tau/\alpha)^{C(\lambda_*+1)}H_\alpha(x,t).
\end{equation}
For $r\le ch$, we take $\alpha=r^2/[A(\lambda_*+1)]$ with $A$ fixed
and sufficiently large. By \eqref{equ-112-8} and \eqref{equ-116-9} at
$\tau=0$,
\[
 \int_{B_{r/2}(x)^c} f(y,t)^2G_\alpha(y-x)\,\mathrm{d}y
 \le\frac{32\alpha}{r^2}\bigl(N_\alpha^x(0)+n/2\bigr)H_\alpha(x,t)
 \le\tfrac12 H_\alpha(x,t).
\]
Thus $H_\alpha(x,t)\le C\alpha^{-n/2}\int_{B_{r/2}(x)}u(y,t)^2\,\mathrm{d}y$.
For $0\le\tau\le4r^2$, the multiplier in \eqref{equ-117-1} is bounded
by $C_*$ and $G_{\alpha+\tau}\ge c_*r^{-n}$ on $B_{2r}(x)$, where
$f=u$. It follows that
\[
 \int_{B_{2r}(x)}u(y,t-\tau)^2\,\mathrm{d}y
 \le C_*r^n H_{\alpha+\tau}(x,t-\tau)
 \le C_*\int_{B_{r/2}(x)}u(y,t)^2\,\mathrm{d}y.
\]
Testing the equation with $\eta\zeta^2u$ gives the terminal energy bound
\[
 \int_{B_{r/2}(x)}u(y,t)^2\,\mathrm{d}y
 \le C(r^{-2}+M)\int_{\Pcal_r(x,t)}u^2
 \le Cr^{-2}\int_{\Pcal_r(x,t)}u^2,
\]
where $\eta(q)=(q-t+r^2)/r^2$ on $[t-r^2,t]$ and $\zeta$ is a cutoff
supported in $B_r(x)$, equal to one on $B_{r/2}(x)$, with
$|\nabla\zeta|\le C/r$.
The constants $C_*$ and $c_*$ depend only on the fixed data and $\lambda_*$.
Integrating in $\tau$, we obtain
\[
 \int_{\Pcal_{2r}(x,t)}u^2\le C_*\int_{\Pcal_r(x,t)}u^2.
\]
Since $\Pcal_r(x,t)\subset\Q_r(x,t)\subset\Pcal_{\sqrt n r}(x,t)$,
a fixed number of iterations gives
\[
 \int_{\Q_{2r}(x,t)}u^2
 \le\int_{\Pcal_{2\sqrt n r}(x,t)}u^2
 \le C_*\int_{\Pcal_r(x,t)}u^2
 \le C_*\int_{\Q_r(x,t)}u^2,
\]
after decreasing $c$. Dividing by the cylinder volumes and choosing
$D_*$ sufficiently large proves \eqref{equ-116-7}.
\end{proof}

\subsection{Propagation of smallness under a doubling condition}\label{subsec-propagation}
Work at terminal time \(0\), and write
\(\Q_r=\Q_r(0,0)\) and \(Q=Q_1(0)\).
Let \(v\) be a real solution of \(v_\tau-\Delta v+cv=0\) in
\(\Q_{8\sqrt n}\), with real \(c\), \(\norm c_\infty\le M_0\), and assume
\begin{equation}\label{equ-117-2}
 \avg_{\Q_{2r}(x,0)}v^2\le4^D\avg_{\Q_r(x,0)}v^2,
 \qquad x\in Q_{4\sqrt n}(0),\quad0<r\le2\sqrt n,
\end{equation}
where \(D\ge0\).
The parabolic degree of \(y^\alpha\tau^j\) is \(|\alpha|+2j\).

\begin{lemma}\label{lem-caloric-polynomials}
Let \(K\ge2\).
If \(p\) has spatial degree at most \(K\) and
\[
 P(y,\tau)=\mathrm{e}^{\tau\Delta}p(y)
   =\sum_{j=0}^{\lfloor K/2\rfloor}\frac{\tau^j}{j!}\Delta^jp(y),
\]
then
\begin{equation}\label{equ-117-3}
 \norm P_{L^2(\Q_1)}
 \le \mathrm{e}^{C_n(K+1)\log(K+2)}\norm p_{L^2(Q)}.
\end{equation}
If \(h_\tau-\Delta h=0\) in \(\Q_{2R}\),
\(R\ge C_n(K+2)^2\), and \(\norm h_{L^\infty(\Q_{2R})}\le B\),
there is such a polynomial \(P\), of parabolic degree at most \(K\), with
\begin{equation}\label{equ-117-4}
 \norm{h-P}_{L^\infty(\Q_1)}
 \le B\left(\frac{C_n(K+2)^2}{R}\right)^{K+1}.
\end{equation}
\end{lemma}
\begin{proof}
We first use the polynomial Markov inequality on the cube to obtain
\(\norm{\partial_i p}_2\le C K^2\norm p_2\).
Indeed, for the normalized Legendre polynomial of degree \(k\), the squared
norm of its derivative is at most \(C(k+1)^3\).
The squared Hilbert--Schmidt norm of differentiation on degrees at most
\(K\) is therefore at most \(C(K+1)^4\).
By taking tensor products, we obtain the same operator bound in each coordinate.
It follows that
\[
 \norm{\Delta^jp}_{L^2(Q)}
 \le(C_nK^4)^j\norm p_{L^2(Q)}.
\]
By integrating the finite expansion of \(P\) over \(-1<\tau\le0\), we obtain
\eqref{equ-117-3}.

To prove \eqref{equ-117-4}, we iterate the interior gradient estimate
through nested cylinders between \(\Q_R\) and \(\Q_{2R}\). This, along with
\(\partial_\tau^jh=\Delta^jh\), gives
\[
 \sup_{\Q_R}|\partial_y^\alpha\partial_\tau^jh|
 \le B\left(\frac{C_nm}{R}\right)^m,
 \qquad m=|\alpha|+2j\ge1.
\]
These estimates hold up to the terminal left trace.
Let \(p\) be the spatial Taylor polynomial of \(h(\cdot,0)\) at the origin
through degree \(K\). We take the Taylor expansion of \(h\) in time at
\(\tau=0\) through order \(J=\lfloor K/2\rfloor\), and then expand each
time derivative \(\partial_\tau^j h(\cdot,0)\) in space through degree \(K-2j\).
The equation identifies the resulting polynomial with \(\mathrm{e}^{\tau\Delta}p\).
Each spatial remainder uses parabolic derivative order \(K+1\), and the
time remainder uses order \(2J+2\in\{K+1,K+2\}\).
By the derivative bounds and the finite Taylor formula, we obtain
\[
 \norm{h-P}_{L^\infty(\Q_1)}
 \le B\,\frac{C_n^{K+1}(K+2)!}{R^{K+1}}.
\]
After increasing \(C_n\), we obtain \eqref{equ-117-4}. This completes the proof.
\end{proof}

The multivariate Remez inequality of Brudnyi and Ganzburg \cite{BG} gives
\[
 \norm p_{L^\infty(Q)}
 \le(C_n/\nu)^K\operatorname*{ess\,sup}_{F}|p|
 \qquad\text{when }|F|\ge\nu|Q|.
\]
Apply it to the subset of \(F\) where
\(|p|\le(2/|F|)^{1/2}\norm p_{L^2(F)}\), whose measure is at least
\(|F|/2\). Thus
\begin{equation}\label{equ-117-5}
 \norm p_{L^2(F)}
 \ge\left(\frac{\nu}{C_n}\right)^{K+1}\norm p_{L^2(Q)},
 \qquad |F|\ge\nu|Q|,\quad\deg p\le K.
\end{equation}

\begin{proposition}[Quantitative propagation of smallness]
\label{prop-propagation}
Under \eqref{equ-117-2}, there is \(C=C(n,M_0)\ge2\) such that,
for every measurable \(F\subset Q\) with \(|F|\ge\nu|Q|\), \(0<\nu\le1\),
\begin{equation}\label{equ-117-6}
 \norm{v(0)}_{L^2(Q)}
 \le \mathrm{e}^{C(D+1)^2\log(D+2)}
       \left(\frac C\nu\right)^{C(D+1)^2}
       \norm{v(0)}_{L^2(F)}.
\end{equation}
\end{proposition}
\begin{proof}
Throughout this proof, all constants depend only on \(n,M_0\).
We put \(A_*=\avg_{\Q_1}v^2\) and may assume \(A_*>0\).
By interior estimates and one doubling step, we have
\begin{equation}\label{equ-117-7}
 \norm{v(0)}_{L^2(Q)}^2\le C4^D A_*.
\end{equation}

\noindent\emph{Step 1. We estimate the local mass and choose the approximation parameters.}

For \(x\in Q\) and dyadic \(0<r\le1\), we iterate \eqref{equ-117-2}
until the radius belongs to \([2\sqrt n,4\sqrt n)\).
The resulting cylinder contains \(\Q_1\), and their volumes are comparable
up to constants depending only on \(n\). Since at most
\(\log_2(1/r)+C_n\) doubling steps are needed, we obtain
\begin{equation}\label{equ-117-8}
 m_x^2:=\avg_{\Q_r(x,0)}v^2
 \ge\mathrm{e}^{-C_n(D+1)}r^{2D}A_*.
\end{equation}
Here we take \(m_x>0\). Every step is admissible, since the center lies in
\(Q\) and each base radius is less than \(2\sqrt n\).

Next, we choose a sufficiently large constant \(A=A(n,M_0)\) and set
\begin{equation}\label{equ-117-9}
 \begin{gathered}
 K=\lceil4(D+1)\rceil,\qquad B=\frac{A(K+2)^2}{\nu},\\
 \eta=B^{-A(K+1)},\qquad R=B^{4A}.
 \end{gathered}
\end{equation}
We take a dyadic \(r\) such that
\[
 \tfrac12 B^{-A^2(K+1)}<r\le B^{-A^2(K+1)}.
\]
Since \(D<(K+1)/4\), a sufficiently large \(A\) ensures
\(R\ge C_n(K+2)^2\) and
\[
 8Rr\le1/2,\qquad 4CM_0r^2R^2(8R)^D\le\eta,
\]
as well as
\begin{equation}\label{equ-118-1}
 2C(8R)^D\left(\frac{C_n(K+2)^2}{R}\right)^{K+1}\le\eta.
\end{equation}
The same choice allows the approximation errors to be absorbed in the
terminal polynomial and Remez estimates below. Moreover,
\begin{equation}\label{equ-118-2}
 \log(1/r)\le
 C(D+1)[\log(D+2)+\log(C/\nu)].
\end{equation}

\noindent\emph{Step 2. We normalize a cube of positive observation density.}

We partition \(Q\) into congruent cubes \(Q_r(x)\).
We choose one, with center \(x_0\), such that
\(|F\cap Q_r(x_0)|\ge\nu|Q_r|\), and let
\[
 U(y,\tau)=m_{x_0}^{-1}v(x_0+ry,r^2\tau).
\]
It satisfies \(\avg_{\Q_1}U^2=1\) and
\[
 U_\tau-\Delta U=-r^2c(x_0+ry,r^2\tau)U.
\]
To bound \(U\) on \(\Q_{2R}\), we take
\(S=2^{\lceil\log_2(4R)\rceil}\in[4R,8R)\).
By the doubling condition and the volume ratio, we have
\[
 \avg_{\Q_{4R}}U^2
 \le2^{n+2}\avg_{\Q_S}U^2
 \le2^{n+2}(8R)^{2D}.
\]
Since \(Sr<8Rr\le1/2\) and \(x_0\in Q\), the corresponding cylinders
satisfy \(\Q_{Sr}(x_0,0)\subset\Q_2\subset\Q_{8\sqrt n}\).
Interior boundedness, with rescaled potential bounded by \(M_0r^2R^2\),
therefore gives
\begin{equation}\label{equ-118-3}
 \norm U_{L^\infty(\Q_{2R})}\le C(8R)^D.
\end{equation}

\noindent\emph{Step 3. We approximate \(U\) by a caloric polynomial.}

Let \(e\) solve
\[
 e_\tau-\Delta e=-r^2c(x_0+ry,r^2\tau)U\quad\hbox{in }\Q_{2R},
\]
with zero lateral and initial data.
By the Dirichlet heat semigroup and the comparison principle, we obtain
\[
 \norm e_{L^\infty(\Q_{2R})}
 \le4M_0r^2R^2\norm U_{L^\infty(\Q_{2R})}\le\eta.
\]
Thus \(h=U-e\) is caloric and
\(\norm h_\infty\le2C(8R)^D\).
By Lemma~\ref{lem-caloric-polynomials} and \eqref{equ-118-1},
there is \(P=\mathrm{e}^{\tau\Delta}p\), \(\deg p\le K\), such that
\[
 \norm{U-P}_{L^\infty(\Q_1)}\le2\eta.
\]
The estimate also holds at \(\tau=0\).
Since \(\avg_{\Q_1}U^2=1\), the polynomial has
\(\norm P_{L^2(\Q_1)}\ge|\Q_1|^{1/2}/2\).
By \eqref{equ-117-3},
\[
 \norm p_{L^2(Q)}
 \ge c_n \mathrm{e}^{-C_n(K+1)\log(K+2)}.
\]

\noindent\emph{Step 4. We apply the Remez inequality and return to the original scale.}

We set \(F'=\{y\in Q:x_0+ry\in F\}\), so \(|F'|\ge\nu|Q|\).
Since \(P(\cdot,0)=p\), the uniform approximation on the terminal face gives
\(\norm{U(0)-p}_{L^2(F')}\le2\eta|F'|^{1/2}\).
By \eqref{equ-117-5} and the choice of \(A\), we obtain
\[
 \begin{split}
 \norm{U(0)}_{L^2(F')}
 &\ge\norm p_{L^2(F')}-2\eta|F'|^{1/2}\\
 &\ge c_n
  \mathrm{e}^{-C_n(K+1)[\log(K+2)+\log(C_n/\nu)]}.
 \end{split}
\]
This, along with \eqref{equ-117-8}, gives
\[
 \begin{split}
 \int_Fv(x,0)^2\,\mathrm{d}x
 &\ge r^n m_{x_0}^2\int_{F'}U(y,0)^2\,\mathrm{d}y\\
 &\ge r^{n+2D}
 \mathrm{e}^{-C(D+1)[\log(D+2)+\log(C/\nu)]}A_* .
 \end{split}
\]
By inserting \eqref{equ-118-2} and
\eqref{equ-117-7} and taking square roots, we obtain
\eqref{equ-117-6}. This completes the proof.
\end{proof}

Parabolic rescaling gives the same estimate on any cube for which the
corresponding enlarged cylinders satisfy the preceding hypotheses, with
the rescaled potential bound. In the construction below,
the doubling exponent is the fixed \(D_*\), so the coefficient in
\eqref{equ-117-6} depends only on the prescribed observation density
and the fixed data.

\subsection{Proof of the observation estimate at one time}\label{subsec-proof-ob-one-time}
\begin{proof}[Proof of Theorem~\ref{thm-ob-one-time}]
We first remove a boundary layer of volume less than \(m/2\) and cover the remaining
compact region by finitely many fixed cubes whose required concentric
enlargements are compactly contained in \(\Omega\). Then there is a cube \(Q_0\) such that
\(|E\cap Q_0|\ge\mu|Q_0|\), where \(\mu=\mu(\Omega,m)>0\).
We apply Proposition~\ref{prop-frequency-covering} to \(E\cap Q_0\).
The resulting cube \(Q_*\) has observation density at least \(\mu/2\),
fixed past frequency \(\lambda_*\), and
\[
 \mathrm{e}^{-C_\mu\sqrt\Lambda}\le R_*\le c/\sqrt\Lambda.
\]
Proposition~\ref{prop-terminal-doubling} supplies fixed terminal cylinder
doubling at scales at most \(h=R_*/[C_\mu(1+\Lambda)^2]\).

Next, we partition \(Q_*\) dyadically into cubes of half-side
\(c_1h\le\ell\le c_2h\).
We choose \(c_2\) small enough that all enlarged cylinders required in
Proposition~\ref{prop-propagation} have centers in \((3/2)Q_*\), radii at most
\(ch\), and time depth at most \(h^2\).
At least one child \(Q\) has
\(|E\cap Q|\ge(\mu/2)|Q|\).
Parabolic rescaling changes the potential to \(\ell^2a\), so
Proposition~\ref{prop-propagation}, with \(D=D_*\) and
\(\nu=\mu/2\), gives
\[
 \norm{u(t)}_{L^2(Q)}
 \le C_\mu\norm{u(t)}_{L^2(E\cap Q)}.
\]
Here one may take
\(C_\mu=\mathrm{e}^{C(D_*+1)^2\log(D_*+2)}(2C/\mu)^{C(D_*+1)^2}\),
which is independent of \(\Lambda\).
Moreover,
\[
 1+|\log\ell|\le C_\mu\sqrt\Lambda+2\log(1+\Lambda)+C_\mu
              \le C_\mu'\sqrt\Lambda.
\]
Thus all three conditions in \eqref{equ-114-1} hold.
By Proposition~\ref{prop-cube-criterion}, we obtain \eqref{equ-111-1}, which completes the proof.
\end{proof}

\section{The whole space: thick slices and integral thickness}\label{sec-whole-space}

\subsection{Whole-space observability}\label{subsec-ob-Rn}\label{subsec-proof-ob-Rn}

In this section we write $\norm{\cdot}_2=\norm{\cdot}_{L^2(\R^n)}$. Put
\[
 K=1000(n+1),\qquad Q=z+[-\ell/2,\ell/2]^n,
 \qquad Q^+=z+[-K\ell,K\ell]^n.
\]
The fixed factor $K$ ensures that the cutoff supports and spatial
enlargements used below lie strictly inside $Q^+$.

\begin{lemma}\label{lem-ob-cube}
Let $0\le s<t\le T$, $\delta=t-s$, and let $u$ be an energy solution of
\eqref{equ-1031-5}.
Assume
\[
 A_Q=\int_Q|u(x,t)|^2\,\mathrm{d}x>0,\qquad
 B_Q=\delta^{-1}\int_s^t\int_{Q^+}|u(x,\tau)|^2\,\mathrm{d}x\,\mathrm{d}\tau.
\]
Define
\[
 \Lambda_Q=1+\delta^{-1}+\log^+(B_Q/A_Q),\qquad
 \log^+r=\max\{0,\log r\}.
\]
For every measurable $E\subset Q$ with $|E|\ge\gamma|Q|$,
\begin{equation}\label{equ-118-4}
 A_Q^{1/2}\le C_\gamma \mathrm{e}^{C_\gamma \Lambda_Q^{3/2}}
                         \norm{u(t)}_{L^2(E)}.
\end{equation}
Here $C_\gamma=C(n,M,T,\gamma,\ell)$.
\end{lemma}
\begin{proof}
\noindent\emph{Step 1. We obtain two local energy estimates.}

We fix $\rho=c_n\min\{1,\ell\}$ with $c_n$ small and take centers near $3Q$.
We choose $R=C_n\ell$ so that $Q\subset B_R(x)$ and all auxiliary spatial
domains lie in $Q^+$. The two-ball--one-cylinder inequality \cite[Lemma 3.2]{DWZ},
after translating time, gives
\begin{equation}\label{equ-118-5}
 A_Q\le[C\mathrm{e}^{C(1+\delta^{-1})}\delta B_Q]^\theta
       \left[C\int_{B_{\rho/4}(x)}|u(y,t)|^2\,\mathrm{d}y\right]^{1-\theta},
 \qquad0<\theta<1.
\end{equation}
Here $\theta$ depends only on the fixed radii. By rearranging this inequality, we obtain
\begin{equation}\label{equ-118-6}
 \int_{B_{\rho/4}(x)}|u(y,t)|^2\,\mathrm{d}y\ge \mathrm{e}^{-C\Lambda_Q}A_Q.
\end{equation}
Interior energy estimates, on a fixed interior subdomain $U\Subset Q^+$ containing all
cutoff supports, give
\begin{equation}\label{equ-118-7}
 \sup_{\tau\in[t-\delta/4,t]}\norm{u(\tau)}_{H^1(U)}^2
       \le C(1+\delta^{-2})\delta B_Q\le \mathrm{e}^{C\Lambda_Q}A_Q.
\end{equation}
Next, we apply the local mass estimate \cite[Lemma 2.3(ii)]{DWZ} with terminal
radius $\rho/2$, historical interval $[t-\delta/2,t]$, and a fixed-size cube
centered at $x$, compactly contained in $Q^+$ and containing a fixed
neighborhood of $Q$. Its historical integral is at
most $\delta B_Q$ and, by an interior energy estimate, at least
$\mathrm{e}^{-C\Lambda_Q}A_Q$. Together with
\eqref{equ-118-6}, this bounds the reciprocal of its time
parameter by $C\Lambda_Q$ and gives
\begin{equation}\label{equ-118-8}
 \int_{B_\rho(x)}|u(y,\tau)|^2\,\mathrm{d}y\ge \mathrm{e}^{-C\Lambda_Q}A_Q,
 \qquad t-c_1/\Lambda_Q\le \tau\le t.
\end{equation}
We take $c_1$ small enough that this interval lies in $[t-\delta/8,t]$.

\noindent\emph{Step 2. We estimate the frequency with a common cutoff.}

We choose $\chi$ supported in $U$, equal to one near $3Q$, and put
$f=\chi u$, $g=(\partial_t-\Delta)f$.
We now use the Gaussian argument of Proposition~\ref{prop-local-frequency} and
Lemma~\ref{lem-common-cutoff} on the time window in \eqref{equ-118-8}.
We choose $c_2$ and the frequency scale constant small enough that all backward
windows based at $\tau\in[t-c_2/\Lambda_Q,t]$ remain in this interval.
It follows from \eqref{equ-118-7} and \eqref{equ-118-8} that
\begin{equation}\label{equ-119-1}
 \begin{aligned}
 &N_x(r)\le C\Lambda_Q,\qquad0<r\le c/\sqrt{\Lambda_Q},\\
 &\frac{\int|g(y,t-r^2)|^2G_{r^2}(y-x)\,\mathrm{d}y}
      {\int|f(y,t-r^2)|^2G_{r^2}(y-x)\,\mathrm{d}y} \le C,
 \end{aligned}
\end{equation}
where $N_x$ is defined by \eqref{equ-113-11} for this cutoff.
The contribution from the potential part of $g$ is bounded by $CM^2$; the cutoff terms are bounded by
$C\exp(C\Lambda_Q-c/r^2)$, which can be absorbed by reducing the scale constant.
The same argument gives \eqref{equ-113-3} with $\Lambda_Q$ at the
nearby times. By iterating from \eqref{equ-118-8}, we obtain, for $0<r\le\rho$,
\begin{equation}\label{equ-119-2}
 \int_{B_r(x)}|u(y,\tau)|^2\,\mathrm{d}y
       \ge \mathrm{e}^{-C\Lambda_Q(1+|\log(r/\ell)|)}A_Q,
 \qquad t-c_2/\Lambda_Q\le \tau\le t.
\end{equation}

\noindent\emph{Step 3. We reduce the frequency and apply propagation of smallness.}

If $Q$ has half-side less than one, take $Q_0=Q$. Otherwise, subdivide $Q$
into congruent dyadic cubes of half-side less than one and choose a cube
$Q_0$ with $|E\cap Q_0|\ge\gamma|Q_0|$. In either case, $3Q_0\subseteq 3Q$.
By \eqref{equ-119-1}, we can apply Proposition~\ref{prop-frequency-covering}
to $E\cap Q_0$ with density $\gamma$.
We thus obtain a cube $Q_*\subset Q$ of half-side $R_*$ such that
\[
 |E\cap Q_*|\ge(\gamma/2)|Q_*|,\qquad
 \mathrm{e}^{-C_\gamma\sqrt{\Lambda_Q}}\le R_*\le c/\sqrt{\Lambda_Q},
\]
and the past frequency on $3Q_*$ is bounded by a fixed $\lambda_*(\gamma)$.
With these frequency bounds, we can argue as in the proof of Proposition~\ref{prop-terminal-doubling},
using \eqref{equ-118-7} and \eqref{equ-119-2}
as its energy inputs; its auxiliary window
$h^2$ in \eqref{equ-116-8} lies in $[0,c_2/\Lambda_Q]$ after increasing its fixed
constant. We therefore obtain terminal doubling with a fixed exponent $D_*$.
We partition $Q_*$ into cubes at a sufficiently small fixed multiple of $h$.
One child $P$, of half-side $r$, has $|E\cap P|\ge(\gamma/2)|P|$ and
$1+|\log(r/\ell)|\le C_\gamma\sqrt{\Lambda_Q}$.
The corresponding enlarged cylinders meet the hypotheses of
Proposition~\ref{prop-propagation}, so
\[
 \norm{u(t)}_{L^2(P)}\le C_\gamma\norm{u(t)}_{L^2(E\cap P)}.
\]
Finally, \eqref{equ-119-2} on the ball inscribed in $P$ gives
\[
 \begin{aligned}
 A_Q^{1/2}&\le C_\gamma
 \mathrm{e}^{C_\gamma\Lambda_Q(1+|\log(r/\ell)|)}\norm{u(t)}_{L^2(P)}\\
 &\le C_\gamma\mathrm{e}^{C_\gamma\Lambda_Q^{3/2}}\norm{u(t)}_{L^2(E)},
 \end{aligned}
\]
which proves \eqref{equ-118-4} and completes the proof.
\end{proof}

\begin{proposition}\label{prop-ob-one-time-Rn}
For every $(\gamma,\ell)$-thick measurable set $\omega$ and every
$0\le s<t\le T$ with $\norm{u(t)}_2>0$, set
\[
 \Theta_{s,t}=\mathrm{e}^{2M(t-s)}\norm{u(s)}_2^2/\norm{u(t)}_2^2,\qquad
 \Lambda=\Lambda_{s,t}.
\]
Then
\begin{equation}\label{equ-119-3}
 \norm{u(t)}_2\le C \mathrm{e}^{C\Lambda^{3/2}}\norm{u(t)}_{L^2(\omega)},
 \qquad C=C(n,M,T,\gamma,\ell).
\end{equation}
\end{proposition}
\begin{proof}
We partition $\R^n$ into half-open cubes $Q_k$ of side $\ell$ and let $Q_k^+$
be the translated enlarged cubes in Lemma~\ref{lem-ob-cube}.
We set
\[
 A_k=\int_{Q_k}|u(x,t)|^2\,\mathrm{d}x,\qquad
 B_k=\frac1{t-s}\int_s^t\int_{Q_k^+}|u(x,\tau)|^2
                    \,\mathrm{d}x\,\mathrm{d}\tau.
\]
The enlarged cubes have overlap at most a fixed integer $N=N(n)$, since
their size relative to $Q_k$ is fixed. For $s\le\tau\le t$, the energy
inequality gives
$\norm{u(\tau)}_2^2\le\mathrm{e}^{2M(t-s)}\norm{u(s)}_2^2$.
Consequently, Tonelli's theorem yields
\begin{equation}\label{equ-119-4}
 \begin{aligned}
 \sum_k A_k&=\norm{u(t)}_2^2,\\
 \sum_k B_k&\le\frac{N}{t-s}\int_s^t\norm{u(\tau)}_2^2\,\mathrm{d}\tau
       \le N\Theta_{s,t}\norm{u(t)}_2^2.
 \end{aligned}
\end{equation}
We define the set of good indices by
\[
 \mathcal G=\{k:A_k>0,\ B_k\le2N\Theta_{s,t}A_k\}.
\]
If $k\notin\mathcal G$ and $A_k>0$, then
$A_k<B_k/(2N\Theta_{s,t})$. Indices with $A_k=0$ contribute no terminal
mass, so \eqref{equ-119-4} implies
\[
 \sum_{k\notin\mathcal G}A_k
 \le\frac1{2N\Theta_{s,t}}\sum_kB_k\le\frac12\norm{u(t)}_2^2.
\]
It follows that the good cubes retain at least half of the terminal energy.
For $k\in\mathcal G$, we use $\Theta_{s,t}\ge1$ to obtain
\[
 \Lambda_{Q_k}=1+(t-s)^{-1}+\log^+(B_k/A_k)
 \le\Lambda+\log(2N)\le C_n\Lambda.
\]
By thickness, we have $|\omega\cap Q_k|\ge\gamma|Q_k|$; cube boundaries have
measure zero. We apply Lemma~\ref{lem-ob-cube} with
$E=\omega\cap Q_k$ and square the resulting inequality to obtain
\[
 A_k\le C\mathrm{e}^{C\Lambda^{3/2}}
       \int_{\omega\cap Q_k}|u(x,t)|^2\,\mathrm{d}x,
 \qquad k\in\mathcal G,
\]
with a constant independent of $k$. Since the cubes are disjoint, we can sum
over the good indices to obtain
\[
 \frac12\norm{u(t)}_2^2\le\sum_{k\in\mathcal G}A_k
 \le C\mathrm{e}^{C\Lambda^{3/2}}\norm{u(t)}_{L^2(\omega)}^2.
\]
By taking square roots and adjusting $C$, we obtain \eqref{equ-119-3}. This ends the proof.
\end{proof}

\begin{proof}[Proof of Theorem~\ref{thm-ob-Rn}]
We apply Proposition~\ref{prop-subquadratic} with
\[
 h(t)=\mathrm{e}^{-Mt}\norm{u(t)}_2,\qquad
 g(t)=\mathrm{e}^{-Mt}\norm{u(t)}_{L^2(D_t)}.
\]
The energy inequality makes $h$ nonincreasing, and
Proposition~\ref{prop-ob-one-time-Rn} supplies its hypothesis uniformly on
$\mathcal A$, with $\beta=3/2$. We thus obtain \eqref{equ-1031-6};
the factor $\mathrm{e}^{MT}$ is absorbed into $C$. This completes the proof.
\end{proof}

\subsection{Failure under integral thickness}\label{subsec-integral-thickness}

Write $S_t=\mathrm{e}^{t\Delta}$. A measurable set
$D\subset\R^n\times(0,T)$ is \emph{integral-thick} if there are $R,\eta>0$ such that
\begin{equation}\label{equ-119-5}
 \int_0^T|D_t\cap B_R(y)|\,\mathrm{d}t\ge\eta
 \qquad\text{for every }y\in\R^n.
\end{equation}
This condition averages the observation geometry over time.

\begin{theorem}[Failure under integral thickness]\label{thm-integral-counterexample}
There is a measurable integral-thick set $D\subset\R\times(0,1)$ for which
no finite $C$ satisfies
\[
 \norm{S_1f}_2^2\le C\int_D|S_tf(x)|^2\,\mathrm{d}x\,\mathrm{d}t
 \qquad\text{for all }f\in L^2(\R).
\]
\end{theorem}
\begin{proof}
\noindent\emph{Step 1. We construct solutions with small early-time values and a nonzero terminal norm.}

We choose a real Schwartz function $\varphi$ with $\norm\varphi_2=1$ and
$\operatorname{supp}\widehat\varphi\subset(-1/4,1/4)$, using
$\widehat{S_tf}(\xi)=\mathrm{e}^{-t\xi^2}\widehat f(\xi)$.
For $m\ge1$, we put
\[
 P_m(z)=(1-z)^m\sum_{j=0}^{m-1}\binom{m+j-1}{j}z^j
       =1+\sum_{j=m}^{2m-1}a_{m,j}z^j.
\]
The last equality follows by truncating the Taylor series of $(1-z)^{-m}$.
Moreover,
\[
 \sum_{j=0}^{m-1}\binom{m+j-1}{j}
       =\binom{2m-1}{m-1}\le4^m.
\]
The sum of the absolute coefficients of $(1-z)^m$ is $2^m$, and
$|1-\mathrm{e}^{i\theta}|\le|\theta|$. Consequently,
\[
 \sum_{j=m}^{2m-1}|a_{m,j}|\le8^m,\qquad
 |P_m(\mathrm{e}^{i\theta})|\le(4|\theta|)^m.
\]
We set $\alpha=128^{-2}$, $N_k=64(k+1)$, $m_k=N_k^2$, and define
\[
 f_k(x)=\varphi(x)\operatorname{Re}P_{m_k}(\mathrm{e}^{4ix/N_k}),
 \qquad h_k(t)=S_tf_k.
\]
For $t>0$, the heat semigroup representation gives
\[
 h_k(t,x)=\operatorname{Re}\int_\R G_t(x-y)\varphi(y)
                    P_{m_k}(\mathrm{e}^{4iy/N_k})\,\mathrm{d}y.
\]
Using the bound on $P_{m_k}$, we obtain
\[
 \begin{aligned}
 |h_k(t,x)|
 &\le\norm\varphi_\infty(16/N_k)^{m_k}
       \int_\R G_t(x-y)|y|^{m_k}\,\mathrm{d}y\\
 &\le\norm\varphi_\infty
       \left[\frac{16}{N_k}(|x|+2\sqrt{tm_k})\right]^{m_k}.
 \end{aligned}
\]
The second inequality follows from Minkowski's inequality and the kernel
moment bound
$\bigl(\int_\R G_t(z)|z|^{m_k}\,\mathrm{d}z\bigr)^{1/m_k}
\le2\sqrt{tm_k}$.
Thus, by continuity also at $t=0$,
\begin{equation}\label{equ-119-6}
 |h_k(t,x)|\le\norm\varphi_\infty
       [16k/N_k+32\sqrt\alpha]^{m_k}
 \le\norm\varphi_\infty2^{-m_k},
 \quad 0\le t\le\alpha,\quad |x|\le k.
\end{equation}
We write $h_k(t)=S_t\varphi+b_k(t)$, where
\[
 b_k(t)=S_t\!\left[\varphi\sum_{j=m_k}^{2m_k-1}
                 a_{m_k,j}\cos(4j\,\cdot/N_k)\right].
\]
Each modulated profile has $L^2$ norm at most one. Since its Fourier support lies
in $\{|\xi|\ge4N_k-1/4\}\subset\{|\xi|\ge3N_k\}$, we obtain
\begin{equation}\label{equ-119-7}
 \norm{b_k(t)}_2\le8^{m_k}\mathrm{e}^{-9tN_k^2}
       \le\mathrm{e}^{-m_k}\qquad(1/2\le t\le1).
\end{equation}
Here the last inequality uses $m_k=N_k^2$ and $\log8-9/2<-1$.
Its Fourier support is disjoint from that of $S_t\varphi$. It follows that
\begin{equation}\label{equ-119-8}
 \norm{h_k(1)}_2^2\ge\norm{S_1\varphi}_2^2\ge\mathrm{e}^{-1/8}.
\end{equation}

\noindent\emph{Step 2. We construct an integral-thick observation set.}

By contractivity, we have
$\int_0^\alpha\norm{h_k(t)}_2^2\,\mathrm{d}t
\le\alpha\norm{f_k}_2^2<\infty$.
The spatial tails therefore tend to zero. For each $k$, we can choose $\rho_k\ge k$
such that
\[
 \int_0^\alpha\int_{|x|>\rho_k}|h_k(t,x)|^2
                         \,\mathrm{d}x\,\mathrm{d}t\le2^{-k}.
\]
We take $c_1=0$, $c_{k+1}=c_k+2k+1+\rho_k+\rho_{k+1}$, and put
\[
 E=\bigcup_{k\ge1}[c_k-k,c_k+k],\qquad
 D=\bigl(E\times(0,\alpha)\bigr)
        \cup\bigl(E^c\times(1/2,1)\bigr).
\]
For every $x\in\R$, the total observation time is
\[
 \int_0^1\mathbf1_D(x,t)\,\mathrm{d}t
 =\alpha\mathbf1_E(x)+\tfrac12\mathbf1_{E^c}(x)\ge\alpha.
\]
By integrating over $(y-1,y+1)$ and using Tonelli's theorem, we obtain
$\int_0^1|D_t\cap(y-1,y+1)|\,\mathrm{d}t\ge2\alpha$ for every $y$.
Thus $D$ is integral-thick.
The gap between consecutive intervals in $E$ is $\rho_k+\rho_{k+1}$.
Consequently, all intervals except the $k$th lie outside
$[c_k-\rho_k,c_k+\rho_k]$, whereas the $k$th is $[c_k-k,c_k+k]$.
Hence
\[
 (E-c_k)\setminus[-k,k]\subset\{|x|>\rho_k\},\qquad
 (E-c_k)^c\subset\{|x|>k\}.
\]
\noindent\emph{Step 3. We show that observability fails.}

At early times, we use \eqref{equ-119-6} on the central interval,
while the remaining intervals contribute at most $2^{-k}$.
At late times, we use $|h_k|^2\le2|S_t\varphi|^2+2|b_k|^2$ and
\eqref{equ-119-7}. Combining these estimates, we obtain
\begin{equation}\label{equ-119-9}
 \begin{aligned}
 \int_D|h_k(t,x-c_k)|^2\,\mathrm{d}x\,\mathrm{d}t
 =\int_0^\alpha\int_{E-c_k}|h_k(t,x)|^2\,\mathrm{d}x\,\mathrm{d}t  +\int_{1/2}^1\int_{(E-c_k)^c}|h_k(t,x)|^2
                                      \,\mathrm{d}x\,\mathrm{d}t\\
 \le2\alpha k\norm\varphi_\infty^2 4^{-m_k}+2^{-k}+\mathrm{e}^{-2m_k}+2\int_{1/2}^1\int_{|x|>k}|S_t\varphi(x)|^2
                         \,\mathrm{d}x\,\mathrm{d}t\longrightarrow0.
 \end{aligned}
\end{equation}
The late-time error satisfies
$2\int_{1/2}^1\norm{b_k(t)}_2^2\,\mathrm{d}t\le\mathrm{e}^{-2m_k}$.
The last integral tends to zero by dominated convergence, since
$\int_{1/2}^1\norm{S_t\varphi}_2^2\,\mathrm{d}t\le1/2$.
The other terms vanish because $m_k=64^2(k+1)^2$.
Finally, we set $p_k(x)=f_k(x-c_k)/\norm{S_1f_k}_2$.
Since the heat semigroup commutes with translations,
$S_tp_k(x)=h_k(t,x-c_k)/\norm{S_1f_k}_2$.
The denominator is at least $\mathrm{e}^{-1/16}$ by
\eqref{equ-119-8}, so normalization increases the observation
integral by at most $\mathrm{e}^{1/8}$. Thus
\eqref{equ-119-9} implies
\begin{equation}\label{equ-119-10}
 \norm{S_1p_k}_2=1,\qquad
 \int_D|S_tp_k(x)|^2\,\mathrm{d}x\,\mathrm{d}t\longrightarrow0,
\end{equation}
whereas the proposed observability estimate would have left-hand side
equal to one and right-hand side tending to zero. This contradiction completes the proof.
\end{proof}

\end{document}